\documentclass[reqno,11pt]{amsart}
\usepackage{amsmath,amsthm,amssymb,amsfonts,amscd}
\usepackage{graphicx,color}
\usepackage{boxedminipage}

\input xy
\xyoption{all}

\newcommand{\R}{{\mathbb{R}}}

\newcommand{\Z}{{\mathbb{Z}}}

\newcommand{\bp}{\begin{pmatrix}}
\newcommand{\ep}{\end{pmatrix}}
\newcommand{\bps}{\begin{smallmatrix}}
\newcommand{\eps}{\end{smallmatrix}}

\def \cC{{\mathcal C}}

\def \cF{{\mathcal F}}
\def \cG{{\mathcal G}}
\def \cH{{\mathcal H}}

\def \cS{{\mathcal S}}

\def \law{\leftarrow}

\def \la{\langle}
\def \ra{\rangle}

\def \grad{\mathrm{grad}}
\def \Id{\mathrm{Id}}

\def\Hom{\mathrm{Hom}}

\def \Ob{\mathrm{Ob}}

\def \Ms{\mathit{Ms}}
\def \ie{{\it i.e.}}

\def \be{{\bf e}}

\def \0{{\bf 0}}
\def \1{{\bf 1}} 

\def \ov#1{\frac{1}{#1}}

\newcommand{\ti}{\tilde}

\def \m{{\frak m}}

\def \gn0#1#2{({#1})_{#2}}
\def \vv{{\overset{\circ}{v}}}
\def \fF{{\frak F}}

\def \Nsddata#1#2#3#4#5{
(\xymatrix{{#1}\  \ar@<0.5ex>[r]^{{#2}} & \ {#4}
\ar@<0.5ex>[l]^{{#3}}} ,#5) }

\numberwithin{equation}{section}

\newtheorem{thm}{Theorem}[section]

\newtheorem{lem}[thm]{Lemma}
\newtheorem{prop}[thm]{Proposition}

\theoremstyle{definition}

\newtheorem{defn}[thm]{Definition}
\newtheorem{exmp}{Example}

\newtheorem{rem}[thm]{Remark}

\newtheorem*{pf}{Proof}

\begin{document}
\hfill RIMS-1585, IHES/M/07/10\\

\title{An $A_\infty$-structure for lines in a plane}

\date{March, 2007}

\dedicatory{Dedicated to Professor Jim Stasheff in honor of his 
70th birthday}

\author{Hiroshige Kajiura}
\address{Research Institute for Mathematical Sciences, 
Kyoto University,
606-8502, Japan}
\email{kajiura@kurims.kyoto-u.ac.jp}

\begin{abstract}
As an explicit example of an $A_\infty$-structure associated to geometry, 
we construct an $A_\infty$-structure for a Fukaya category 
of finitely many lines (Lagrangians) in $\R^2$, \ie, we define 
also {\em non-transversal} $A_\infty$-products. 
This construction is motivated by homological mirror symmetry of (two-)tori, 
where $\R^2$ is the covering space of a two-torus. 
The strategy is based on an algebraic reformulation of Morse homotopy 
theory through homological perturbation theory (HPT) 
as discussed by Kontsevich and Soibelman in \cite{KoSo}, 
where we introduce a special DG category which is a key idea 
of our construction. 
\end{abstract}

\maketitle
%
%

\tableofcontents

 \section{Introduction}

\begin{figure}[h]
 \includegraphics{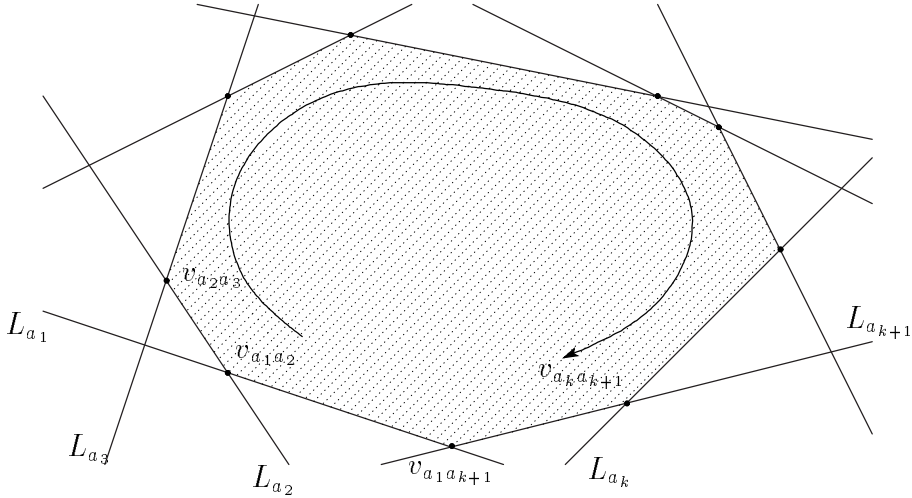}
 \caption{A clockwise convex polygon (CC-polygon) 
defined by lines $L_{a_1},\dots,L_{a_{k+1}}$. }\label{fig:polygon-intro}
\end{figure}

For a graded vector space $A$, a strong homotopy associative structure 
(or an $A_\infty$-structure) on $A$ 
is a family of multilinear maps $m_k:A^{\otimes k}\to A$ for $k\ge 1$ 
satisfying certain constraints, first introduced by Jim Stasheff 
\cite{jds:hahI,jds:hahII} 
in the study of H-spaces such as based loop spaces. 
In particular, $m_1=d$ forms a differential on $A$, $m_2$ is a 
product which is associative up to homotopy, where $m_3$ defines the homotopy 
and $m_4,m_5,\dots$ define higher homotopies. 
An $A_\infty$-algebra with higher products $m_3,m_4,\dots$ all zero is 
a differential graded (DG) algebra, which appears as the structure
DeRham complexes have in general. 
A category version of an $A_\infty$-algebra $(A,\{m_k\}_{k\ge 1})$ 
is called an $A_\infty$-category introduced by Fukaya \cite{fukayaAinfty} 
to formulate Morse homotopy theory and 
Floer theory of Lagrangian submanifolds (shortly, Lagrangians) 
in a symplectic manifold. 
In particular, a category for the latter theory, called a Fukaya category, 
can give an interesting example 
of $A_\infty$-structures associated to geometry. 
In this paper, we shall construct an $A_\infty$-structure 
for a Fukaya category $Fuk(\R^2)$ consisting of lines in $\R^2$. 
For each two lines $L_a$, $L_b$ which intersect with each other 
at one point $v_{ab}\in\R^2$, 
the space of morphisms $\Hom(a,b)$ is a one-dimensional vector space 
(over $\R$) 
spanned by a base $[v_{ab}]$ associated to the intersection point $v_{ab}$. 
Then, the (higher) $A_\infty$-product 
$m_k:\Hom({a_1},{a_2})\otimes\cdots\otimes\Hom({a_k},{a_{k+1}})
\to\Hom({a_1},{a_{k+1}})$, $a_1,\dots,a_{k+1}\in\Ob(\cC)$, 
is defined by polygons surrounded by lines in $\R^2$ 
\begin{equation*}
 m_k([v_{a_1a_2}],\dots,[v_{a_ka_{k+1}}])
 =\pm e^{-Area(\vec{v})} [v_{a_1a_{k+1}}]
\end{equation*}
if the sequence 
$\vec{v}:=(v_{a_1a_2},\dots.,v_{a_ka_{k+1}},v_{a_{k+1}a_1})$, 
$v_{a_{k+1}a_1}=v_{a_1a_{k+1}}$, 
of the intersection points forms 
a clockwise convex (CC-) polygon (Figure \ref{fig:polygon-intro}). 
From the viewpoint of Lagrangian intersection Floer theory, 
these lines are thought of as special Lagrangian submanifolds 
in a symplectic manifold $\R^2\simeq T^*\R$, the cotangent bundle over $\R$. 
We in particular construct such an $A_\infty$-category $Fuk(\R^2)$ 
with finitely many objects, 
in which $\oplus_{a,b\in\Ob(Fuk(\R^2))}\Hom(a,b)$ is 
an example of an $A_\infty$-algebra. 
Although the definition of an $A_\infty$-structure of a Fukaya category 
is clear for the multilinear maps $m_k$ 
on morphisms $\Hom(a,b)$ with $L_a$ and $L_b$ transversal to each other, 
even in this $\R^2$ case, 
it is technically not easy to define 
multilinear maps $m_k$ on morphisms including 
$\Hom(a,a)$ for some line $L_a$ because of non-transversality 
of Lagrangians. (See FOOO \cite{FOOO} 
for the problem of transversality in a more general setup, where, 
I have heard, another way of resolution is discussed.)
However, we can not define an $A_\infty$-category 
without defining all those non-transversal multilinear maps. 
The aim of this paper is to define explicitly all the 
$A_\infty$-products of the Fukaya category $Fuk(\R^2)$ 
including these non-transversal ones. 
To derive those $A_\infty$-products, 
the rough direction of our strategy is 
first to define a DG category $\cC_{DR}$ with the same objects, 
and then to apply to $\cC_{DR}$ 
homological perturbation theory (HPT) 
developed by 
Gugenheim, Lambe, Stasheff, Huebschmann, Kadeishvili, etc., 
\cite{gugjds, GLS:chen, GLS, HK} 
(see also the decomposition theorem in \cite{hk:thesis, KaTe}). 
For a DG algebra or an $A_\infty$-algebra $A$, 
HPT starts with what is called {\em strong deformation retract (SDR)} data 
\begin{equation*}
\Nsddata{B}{\ \iota}{\ \pi}{A}{h}, 
\end{equation*}
where $B$ is a complex, $\iota$ and $\pi$ are chain maps so that 
$\pi\circ\iota=\Id_B$ and 
$h:A\to A$ is the contracting homotopy defined by 
\begin{equation*}
 d_Ah+hd_A=\Id_A-P, \qquad P:=\iota\circ\pi, 
\end{equation*}
together with some additional conditions. 
By definition, $P$ will be an idempotent in $A$. 
Given a contracting homotopy $h$, one obtains SDR data as above, and 
then the HPT machinery gives a way to produce an $A_\infty$-structure on $B$ 
which is homotopy equivalent to the original $A_\infty$-algebra $A$. 
In particular, the induced $A_\infty$-structure on $B$ can be described 
in terms of planar rooted trees (Feynman graphs). 
The HPT for $A_\infty$-algebras is extended straightforwardly 
to $A_\infty$-categories \cite{KoSo}. 
Since the $A_\infty$-structure is manifest on $\cC_{DR}$ 
as its DG-category structure 
and there are not subtleties about transversality there, 
one can expect that HPT yields an $A_\infty$-structure on $Fuk(\R^2)$ 
if we can find a suitable contracting homotopy $h$. 

Physically, this DG category $\cC_{DR}$ 
is related to a kind of Chern-Simons field theory 
(on one-dimensional space $\R$). 
Applying HPT to $\cC_{DR}$ then corresponds to 
considering perturbation theory of the Chern-Simons theory at tree level. 
This kind of Chern-Simons theory is thought of as a topological 
open string field theory (SFT) \cite{Wi:CS}, where 
the choice of a homotopy operator in applying HPT corresponds to 
the choice of a gauge fixing for the open SFT 
(see \cite{hk:SFT,hk:thesis} for open SFT and 
\cite{lazaroiu:sft-brane, costello:CS} for topological open SFT ). 
From such a physical viewpoint, it is interesting that 
the result of this paper indicates 
the (disk) instantons, which are nonperturbative effects in string theory, 
are also derived by perturbation theory of string {\em field} theory.  

\begin{figure}[h]
  \begin{minipage}[c]{75mm}{
\begin{center}

\vspace*{0.5cm}

 \includegraphics{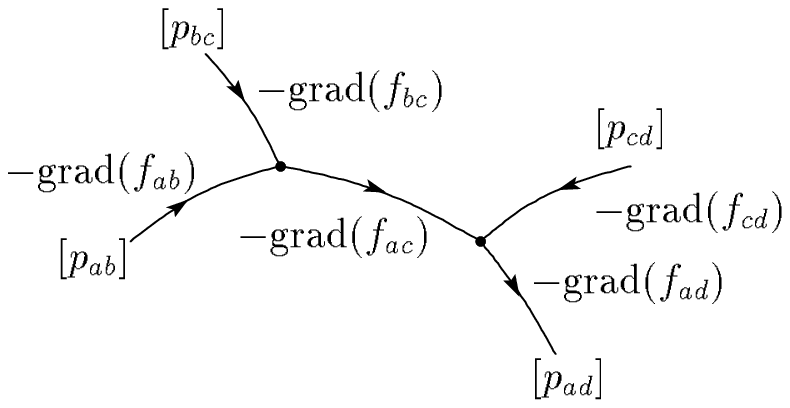}

\vspace*{0.5cm}

(a)
\end{center}}
 \end{minipage}
\quad
 \begin{minipage}[c]{75mm}{
\begin{center}
 \includegraphics{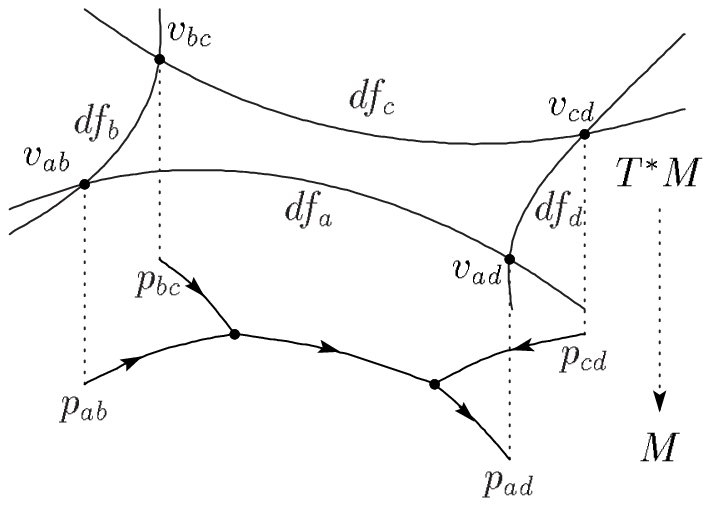}

 (b)
\end{center}}
 \end{minipage}
 \caption{(a): 
A tree of gradient flows in $M$, 
where $f_{ab}:=f_a-f_b$, while $p_{ab}$ is a critical point of $f_{ab}$ 
and $[p_{ab}]$ the associated base of $\Hom_{\Ms(M)}(a,b)$. 
An $A_\infty$-product $m_3([p_{ab}],[p_{bc}],[p_{cd}])$ is defined by 
counting all such trees of gradient flows. \ 
(b): The Lagrangians in $T^*M$ defined by $df_a, df_b, df_c, df_d$ 
corresponding to the Morse functions $f_a, f_b, f_c, f_d$. 
By definition, the dimension of the Lagrangians is the same as 
the dimension of $M$. 
By the projection $x: T^*M\to M$, 
an intersection point $v_{ab}$ of $L_a$ with $L_b$ corresponds to 
a critical point $p_{ab}=x(v_{ab})$. 
}\label{fig:FO}
\end{figure}

The homotopy equivalence $Fuk(\R^2)\simeq \cC_{DR}$ 
obtained via the HPT 
plays the key role in discussing 
homological mirror symmetry \cite{HMS}, 
since the DG category $\cC_{DR}$ is related to a category of holomorphic 
vector bundles on a complex manifold. 
Since $\R^2$ is the covering space of a two-torus, 
the arguments in this paper are directly applied to 
homological mirror symmetry for two-tori, and 
higher dimensional generalization 
of the torus analog of the DG-category $\cC_{DR}$ 
is also straightforward (for instance see \cite{hk:nctheta}). 
Homological mirror symmetry is discussed positively 
for two-tori \cite{PZ, Poli:Ainf}, for abelian varieties \cite{Fabelian}, 
and for (complex) noncommutative tori 
\cite{hk:foliation, PoSc, hk:nchms, hk:ncDG, hk:nctheta}; 
in particular, for two tori, 
transversal $A_\infty$-products are defined explicitly and 
the homological mirror is also shown for the transversal $A_\infty$-products 
by Polishchuk \cite{Poli:Ainf}. 
However, the reason why such equivalence holds has still been unclear 
even for the transversal $A_\infty$-products. 
Kontsevich-Soibelman \cite{KoSo} then proposed 
a strategy to show the homological mirror symmetry based on the viewpoint 
of Strominger-Yau-Zaslow torus fibrations 
(see also \cite{fukaya:asymptotic} for a related approach). 
The strategy is to reformulate 
Fukaya-Oh Morse homotopy theory \cite{fukayaAinfty,FO} 
algebraically in terms of a DG category $DR(M)$ 
consisting of DeRham complexes and to apply HPT to $DR(M)$ 
together with Harvey-Lawson's Morse theory \cite{HaLa}. 
For a compact manifold $M$ with a given metric 
(which is used to define the gradient $\grad$, see below), 
the objects of the category $\Ms(M)$ 
\footnote{This category $\Ms(M)$ is denoted by ${\mathrm M}(Y)$ 
in \cite{KoSo} where 
$Y$ is the compact smooth manifold $M$ here. }
of Fukaya-Oh Morse homotopy \cite{fukayaAinfty,FO} are 
smooth functions $f\in C^\infty(M)$ on $M$. 
If the difference $f_{ab}:=f_a-f_b$ of two functions $f_a, f_b\in C^\infty(M)$ 
is a Morse function, the space $\Hom_{\Ms(M)}(a,b)$ of morphisms is defined as 
the vector space spanned by bases $[p_{ab}]$ 
associated to the critical points $p_{ab}$ of $f_{ab}$. 
The $A_\infty$-structure on $\Ms(M)$ 
is defined by trivalent planar trees so that 
each edge is associated to the gradient flow of the difference of 
the corresponding two functions (see Figure \ref{fig:FO} (a)). 
The equivalence of the Morse $A_\infty$-category $\Ms(M)$ 
with the Fukaya $A_\infty$-category $Fuk(T^*M)$ on $T^*M$ 
is discussed in \cite{FO}, 
where an object of $Fuk(T^*M)$ is a Lagrangian $L_a\subset T^*M$ 
defined by the graph of $df\in\Gamma(T^*M)$ of a Morse function $f_a$ 
(Figure \ref{fig:FO} (b)), and 
the space $\Hom_{Fuk(T^*M)}(a,b)$ of morphisms from $L_a$ to $L_b$ is 
spanned by the bases $[v_{ab}]$ 
associated to the intersection points $v_{ab}$ of $L_a$ with $L_b$, 
whose images by the projection $x: T^*M\to M$ 
are the critical points $p_{ab}=x(v_{ab})$ of $f_{ab}=f_a-f_b$. 
Kontsevich-Soibelman \cite{KoSo} 
discussed obtaining the Morse $A_\infty$-category $\Ms(M)$ 
by applying HPT to the DG-category $DR(M)$. 
The key idea there is to identify 
the contracting homotopy of the SDR in HPT 
with Harvey-Lawson's chain homotopy in \cite{HaLa}, 
which will allow us to identify the planar trees in HPT 
with the trees of gradient flows defining the $A_\infty$-structure of $\Ms(M)$. 
Let $\varphi_t:M\to M$, $t\in [0,\infty)$ be the flow 
generated by the gradient $\grad(f)$ of a given Morse function $f$. 
Harvey-Lawson \cite{HaLa} showed the existence of the limit 
${\bf P}:=\lim_{t\to\infty}\varphi_t^*$ of the pullback $\varphi_t^*$ 
together with the chain homotopy 
\begin{equation}\label{HaLa}
 d_{D'}h+ h d_{DR}={\bf I} -{\bf P}, 
\end{equation}
where ${\bf I}: (\Omega(M),d_{DR})\to (D'(M),d_{DR'})$ is the inclusion of 
smooth differential forms on $M$ to the space $D'(M)$ of distribution forms 
and ${\bf P}:(\Omega(M),d_{DR})\to (D'(M),d_{D'})$ turns out to be a linear map 
such that ${\bf P}(\Omega(M))\subset D'(M)$ forms 
a subcomplex spanned by DeRham currents $[U_{p}]$ with support 
the unstable manifolds $U_{p}$ of critical points $p$ of $f$. 
In \cite{KoSo}, 
the SDR data for the complex $\Hom_{DR(M)}(a,b):=\Omega(M)$ was identified with 
the Harvey-Lawson's chain homotopy (\ref{HaLa}) with $f=f_{ab}:=f_a-f_b$. 
All these tools for $T^*M$ were then extended 
to torus fibrations over $M$ 
to discuss homological mirror symmetry for torus fibrations.

Strongly motivated by this story, 
we define a DG-category $\cC_{DR}$ which is 
similar to $DR(M)$ in \cite{KoSo} with $M=\R$.
Here, for two Morse functions $f_a,f_b\in C^\infty(\R)$, 
we set the differential 
$d_{ab}:\Hom_{\cC_{DR}}(a,b)\to \Hom_{\cC_{DR}}(a,b)$ 
as the twisted differential 
\begin{equation*}
 d_{ab}=d-d(f_{ab})\wedge\ =e^{f_{ab}}\, d\, e^{-f_{ab}}
\end{equation*}
of Witten's Morse complex \cite{Wi:morse}. This leads to the correct 
structure constant of the transversal $A_\infty$-products, 
\ie, the area of the corresponding CC-polygons, 
via the HPT. 
Though the case $M=\R$ looks too simple, 
because $\R$ is noncompact, this case can include 
more nontrivial phenomena than the traditional setting where 
$M$ is a {\em compact} smooth manifold. 
We consider a set $\fF_N:=\{f_a,f_b,\dots\}=\{a,b,\dots\}$ 
of $N$ lines, and denote by $\cC_{DR}(\fF_N)$ 
the DG-category $\cC_{DR}$ with $\Ob(\cC_{DR})=\fF_N$. 
As mentioned above, this theorem is motivated by the case 
where $\R^2\simeq T^*\R$ is replaced by a two-torus 
and its higher dimensional generalizations. 
Our choice of this DG category $\cC_{DR}(\fF_N)$ then comes from 
the DG category of holomorphic vector bundles on a noncommutative 
torus with the noncommutativity set to be zero, 
which gives (an equivalent but) different description from 
the usual commutative torus setting (cf. \cite{PZ, Poli:Ainf}). 
For our purposes, this noncommutative tori setting fits better 
even in discussing commutative tori. 
In particular, we identify $\Omega(M)$, $M=\R$, 
with the space of {\em rapidly decreasing} smooth differential forms; 
for instance, $\Omega^0(\R)$ is the space $\cS(\R)$ of Schwartz functions 
instead of $C^\infty(M)$. 

The main theorem (Theorem \ref{thm:main}) of this paper 
is to show 
the existence of an $A_\infty$-structure on the category $Fuk(\R^2,\fF_N)$
of lines in $\R^2$ whose transversal $A_\infty$-products are 
those associated with CC-polygons 
as in Figure \ref{fig:polygon-intro} 
and which is homotopy equivalent to the DG-category $\cC_{DR}(\fF_N)$. 
We also present explicitly 
one such $A_\infty$-category, which we denote $\cC(\fF_N)$. 
As stated previously, the rough idea to obtain $\cC(\fF_N)$ 
is to apply HPT to $\cC_{DR}(\fF_N)$. 
In order to reproduce the area of the CC-polygon as 
a structure constant of an $A_\infty$-product of $\cC(\fF_N)$, 
the contracting homotopy $h$ of SDR in HPT should be of the type in 
the identity (\ref{HaLa}). 
However, unfortunately, the $h$ in the identity (\ref{HaLa}) 
is the chain homotopy between ${\bf I}$ and ${\bf P}$, which map 
$\Omega(M)$ to not $\Omega(M)$ itself but to $D'(M)$, 
since the DeRham currents $[U_p]$ are {\em not smooth} differential forms. 
Thus, we need some modification of the story. 
One natural way may be to modify $h$ as $h_{\epsilon}$ with a parameter $\epsilon$ 
such that $d_{DR}h_\epsilon+h_\epsilon d_{DR}=\Id-P_\epsilon$ holds on $\Omega(M)$ 
if $\epsilon\ne 0$ and $\lim_{\epsilon\to 0}h_{\epsilon}=h$. 
Then, we may apply HPT with contracting homotopy $h_\epsilon$, $\epsilon\ne 0$, 
construct the induced $A_\infty$-products, 
and finally take the limit $\epsilon\to 0$. 

One such modification $h_\epsilon$ is discussed in \cite{KoSo},
but the strategy in the present paper 
is instead to define a suitable subcomplex of $D'(\R)$. 
Though $D'(\R)$ can not be equipped with a product structure, 
we can introduce a product structure
\footnote{The product structure we shall introduce is also 
motivated by one such modification $h_\epsilon$ which is however different from 
the one discussed in \cite{KoSo}. We hope to discuss the limit $\epsilon\to 0$ 
in this approach elswhere. }
 in the subcomplex 
and apply HPT directly to the subcomplex. 
More precisely, we introduce a DG-category $\cC'_{DR}(\fF_N)$ 
as the smallest DG-category with the same objects 
$\Ob(\cC_{DR}'(\fF_N))=\Ob(\cC_{DR}(\fF_N))=\fF_N$ so that, 
for any $a\ne b\in\fF_N$, 
$\Hom_{\cC_{DR}'(\fF_N)}(a,b)$ includes $[U_{p_{ab}}]\subset D'(\R)$ 
for any critical point $p_{ab}=x(v_{ab})$ of $f_{ab}=f_a-f_b$ 
and is closed with respect to the operation $h$. 
Note that the latter requirement 
enables us to apply HPT directly to $\cC'_{DR}(\fF_N)$. 
This $\cC'_{DR}(\fF_N)$ is in fact homotopy equivalent 
to the original DG-category $\cC_{DR}(\fF_N)$ 
as shown in subsetion \ref{ssec:proof}. 

In our case, $M=\R$ and the graph $df_a$ of $f_a$ is a line $L_a$ 
in $\R^2\simeq T^*\R$ for any $a\in\fF_N$.  
Thus, for any $a\ne b\in\fF_N$, 
the intersection point $v_{ab}$ of $L_a$ and $L_b$ is only one and 
so is the critical point $p_{ab}=x(v_{ab})$ of $f_{ab}$. 
Then, 
$[U_{p_{ab}}]$ will be the gaussian 
$e^{f_{ab}}\in\Omega^0(\R)\subset {D'}^0(\R)$ 
whose support is $\R$ itself 
(but multiplied by $e^{f_{ab}}$ due to 
our choice of the differential $d_{ab}$) 
if the Hessian of $-f_{ab}$ is positive 
and the delta function one-form $\delta_{v_{ab}}\in {D'}^1(\R)$ 
with support $p_{ab}$ if the Hessian of $-f_{ab}$ is negative. 
In order for $\cC'_{DR}(\fF_N)$ to be closed with respect 
to the composition of morphisms, for any $a\ne b\in\fF_N$ 
we need to include $\delta_{v_{cd}}\in\Hom_{\cC'_{DR}(\fF_N)}(a,b)$ for 
any $c\ne d\in\fF_N$. The operation of $h$ on $\delta_{v_{cd}}$ 
will then produce step functions $\vartheta_{v_{cd}}$. 
Consequently, it turns out that the DG category 
$\cC'_{DR}(\fF_N)$ is generated 
by step functions and delta function one forms. 
The contracting homotopy $h_{ab}$ of the type in (\ref{HaLa}) 
gives a desirable idempotent 
$P_{ab}:\Hom_{\cC'_{DR}(\fF_N)}(a,b)\to\Hom_{\cC'_{DR}(\fF_N)}(a,b)$ 
such that $P_{ab}\Hom_{\cC'_{DR}(\fF_N)}(a,b)
=\R\cdot [U_{p_{ab}}]\simeq\R\cdot [v_{ab}]$ 
for $a\ne b\in\fF_N$. 
Here, the corresponding SDR gives a Hodge decomposition of 
$\Hom_{\cC'_{DR}(\fF_N)}(a,b)$. 
However, as we will mention also in the final section, 
there is no natural choice of the Hodge decomposition, 
\ie, contracting homotopy $h_{ab}$ if $a=b$. 
Therefore, we set $h_{aa}=0$. 
Consequently, the space $\Hom_{\cC(\fF_N)}(a,a)$ will be 
a commutative DG algebra (denoted by $A_S(\R)$) which is also 
generated by the step functions and the delta function one-forms. 

This paper is organized as follows. 
After recalling terminologies for $A_\infty$-categories 
and HPT in section \ref{sec:recall}, 
we present the main theorem (Theorem \ref{thm:main}) 
in subsection \ref{ssec:main}. 
Before proving it in section \ref{sec:proof}, 
we present the $A_\infty$-category $\cC(\fF_N)$ explicitly 
in subsection \ref{ssec:explicit}. 
To define the $A_\infty$-category $\cC(\fF_N)$, 
we introduce the commutative DG algebra $A_S(\R)$, which is prepared 
in subsection \ref{ssec:CDGA}. 
Section \ref{sec:interpret} is devoted to presenting 
geometric interpretations of some basic properties 
of the transversal part of the Fukaya $A_\infty$-category $\cC(\fF_N)$ 
in some examples. 
Thus, the contents in section \ref{sec:interpret} 
may essentially be known to experts. 
In subsection \ref{ssec:polygon-degree}, 
we observe that a transversal $A_\infty$-product can be nonzero 
if and only if the corresponding lines form a CC-polygon 
as in Figure \ref{fig:polygon-intro}. 
In subsection \ref{ssec:polygon-relation}, we see 
an $A_\infty$-constraint for transversal $A_\infty$-products 
consists of only two terms which correspond to the ways to 
divide a clockwise polygon with one nonconvex vertex into two. 
We also include a reason why we can not avoid 
non-transversal $A_\infty$-products in subsection \ref{ssec:why}. 
Then, in section \ref{sec:proof} 
we prove the main theorem (Theorem \ref{thm:main}). 
In subsection \ref{ssec:proof}, 
we introduce the DG-category $\cC'_{DR}(\fF_N)$ and 
prove Theorem \ref{thm:main} 
assuming a proposition (Proposition \ref{prop:key}). 
Then, in subsection \ref{ssec:derive} 
we prove Proposition \ref{prop:key}, 
where we derive the $A_\infty$-category $\cC(\fF_N)$ 
by applying HPT to $\cC'_{DR}(\fF_N)$. 
Several examples of the explicit calculations of the derived 
$A_\infty$-products are also given there. 
Since we consider the case $M=\R$, the trees of gradient flows in $M$ 
in the sense of $\Ms(M)$ are degenerate 
to be intervals and points on them. 
On the other hand, the HPT suggests the use of planar trees 
which are useful to determine the signs of the $A_\infty$-products, too. 
Thus, in those calculations, we introduce planar trees 
associated to CC-polygons 
which are lifts of the trees of gradient flows in $M=\R$ to $T^*\R$. 
Finally, we end with mentioning applications of the main theorem 
to the case of tori, etc., in section \ref{sec:conclude}. 

Throughout this paper, 
by (graded) vector spaces 
we indicate those over fields $k=\R$. 
Though motivated strongly by the background stated above, 
the body of this paper can be read independently. 

\noindent
{\bf Acknowledgments} :\ 
First of all, I would like to thank Jim Stasheff 
for his continuous encouragement and valuable 
discussions. I remember communications with him started when 
I discussed some application of HPT to open string field theory; 
the present work is thought of a topological string analog. 
As for discussions related to the background of this work, 
I would also like to thank 
M.~Akaho, T.~Kondo and Y.~Terashima. 
K.~Saito and A.~Takahashi pointed out many faults 
in the formulation of an earlier stage, 
which helped to arrive at the present formulation. 
I am grateful to K.~Fukaya who 
called my attention to various issues about transversality. 
It is needless to say that the present work is motivated 
greatly by what I have learnt from him. 
This work was completed during my visit in IHES 
where the environment for research is excellent and 
I would like to thank all researchers and staff there.


 \section{$A_\infty$-categories and their homotopical properties}
\label{sec:recall}

\subsection{$A_\infty$-algebras}
\label{ssec:Ainfty}

\begin{defn}[$A_\infty$-algebra (strong homotopy associative algebra) 
\cite{jds:hahI,jds:hahII}]
\label{defn:Ainfty}
An $A_\infty$-algebra $(V,\m)$ consists of a 
$\Z$-graded vector space $V$ 
with a collection of 
multilinear maps $\m:=\{m_n:V^{\otimes n}\to V\}_{n\ge 1}$ of 
degree $(2-n)$ satisfying 
\begin{equation}\label{Ainfty}
 \begin{split}
 & 0=\sum_{k+l=n+1}\sum_{j=0}^{k-1}
(-1)^\sigma  
 \ m_k(w_1,\dots,w_j, m_l(w_{j+1},\dots,w_{j+l})
 ,w_{j+l+1},\dots,w_n)\ ,\qquad n\ge 1
 \end{split}
\end{equation}
for homogeneous elements $w_i\in V$, $i=1,\dots,n$, with 
degree $|w_i|\in\Z$, 
where $\sigma=(j+1)(l+1)+l(|w_1|+\cdots+|w_j|)$. 
\end{defn}
That the multilinear map $m_k$ has degree $(2-k)$ indicates 
the degree of $m_k(w_1,\dots,w_k)$ is $|w_1|+\cdots +|w_k|+(2-k)$. 

For $m_1=d$, $m_2=\cdot$, 
the first three relations of the above $A_\infty$-condition are:
\begin{align*}
& d^2=0\ ,\\
& d(w\cdot w')= d(w)\cdot w' + (-1)^{|w|}w\cdot d(w')\ ,\\
& (w\cdot w')\cdot w''-w\cdot (w'\cdot w'')=
 d(m_3)(w,w',w''), \\
 & \qquad\quad d(m_3):=
 d m_3+m_3 (d\otimes\1\otimes\1+\1\otimes d\otimes\1+\1\otimes\1\otimes d)
\end{align*}
for homogeneous elements $w,w',w''\in V$. 
The first identity implies that $(V, d)$ defines a complex. 
The second identity implies that the differential 
$d$ satisfies the Leibniz rule with respect to the product $\cdot$. 

The third identity implies that the product $\cdot$ 
is associative up to homotopy. 
In particular, 
the product $\cdot$ is strictly associative if $m_3$=0. 
\begin{defn}\label{defn:DGA}
An $A_\infty$-algebra $(V,\m)$ with vanishing higher products 
$m_3=m_4=\cdots =0$ is called 
a {\em differential graded algebra (DGA)}. 
\end{defn}

There exists a different definition of $A_\infty$-algebras 
via a shift in degree. 
\begin{defn}\label{defn:Ainfty2}
An $A_\infty$-algebra $(\cH,\m)$ consists of a 
$\Z$-graded vector space $\cH$ 
with a collection of {\em degree one} 
multilinear maps 
$\m:=\{m_n:\cH^{\otimes n}\to \cH\}_{n\ge 1}$ 
satisfying 
\begin{equation*}
 \begin{split}
 & 0=\sum_{k+l=n+1}\sum_{j=0}^{k-1}
(-1)^{|o_1|+\cdots +|o_j|} 
 \ m_k(o_1,\dots,o_j, m_l(o_{j+1},\dots,o_{j+l})
,o_{j+l+1},\dots,o_n)\ .
 \end{split}
\end{equation*}
\end{defn}
These two definitions of $A_\infty$-algebras are in fact equivalent. 
They are related by a degree shifting operator 
\begin{equation*}
 s: V^r\to (V[1])^{r-1}=:\cH^{r-1}\ 
\end{equation*}
called the {\em suspension}. 
The direct relation between multilinear maps in these two definitions is 
given \cite{gj:cyclic} by 
\begin{equation*}
 m^\cH_n=(-1)^{\sum_{i=1}^{n-1}(n-i)}s m^V_n((s^{-1})^{\otimes n})
\end{equation*} 
or more explicitly: 
\begin{equation}\label{suspension-sign}
m^\cH_n(o_1,\dots,o_n) =(-1)^{\sum_{i=1}^{n-1}(n-i)|o_i|}
\ s m^V_n(s^{-1}(o_1),\dots,s^{-1}(o_n)), 
\end{equation}
where we denoted the multilinear maps of $(V,\m)$ and that of $(\cH,\m)$ 
by $\m^V$ and $\m^\cH$, respectively. 

The original definition in Definition \ref{defn:Ainfty} 
is natural in the sense that the differential $m_1$ has degree one, 
the product $m_2$ preserves the degree and then 
$m_n$, $n\ge 3$, are the higher homotopies. 
However, one can see that Definition \ref{defn:Ainfty2} is simpler in sign. 
\begin{defn}[$A_\infty$-morphism]
Given two $A_\infty$-algebras $(\cH,\m)$ and $(\cH',\m')$,  
a collection of degree preserving (= degree zero) 
multilinear maps $\cG:=\{g_k:\cH^{\otimes k}\to\cH'\}_{k\ge 1}$, 
is called an $A_\infty$-morphism $\cG :(\cH,\m)\to (\cH',\m')$ 
if and only if the following relations hold: 
\begin{equation}\label{Ainftymorp}
\sum_{i}\sum_{k_1+\cdots +k_n=n}
 m_i(g_{k_1}\otimes\cdots\otimes g_{k_i}) 
=\sum_{i+1+j=k}\sum_{i+l+j=n} 
g_k(\1^{\otimes i}\otimes m_l\otimes\1^{\otimes j}) 
\end{equation}
for $n=1,2,\dots$. 
\end{defn}
The above relation for $n=1$ implies that 
$g_1:\cH\to\cH'$ forms a chain map $g_1:(\cH,m_1)\to (\cH',m'_1)$. 
\begin{defn}
An $A_\infty$-morphism $\cG:(\cH,\m)\to (\cH',\m')$ is 
called an {\em $A_\infty$-quasi-isomorphism}  
if and only if $g_1:(\cH,m_1)\to (\cH',m_1')$ induces an isomorphism between 
the cohomologies of these two complexes. 
In this situation, we say {\em $(\cH,\m)$ is homotopy equivalent to $(\cH',\m')$} 
and call the $A_\infty$-quasi-isomorphism $\cG:(\cH,\m)\to (\cH',\m')$ 
{\em homotopy equivalence}. 
\end{defn}
It is known that 
there exists an inverse $A_\infty$-quasi-isomorphism 
$\cG':(\cH',\m')\to (\cH,\m)$ for a given 
$A_\infty$-quasi-isomorphism $\cG:(\cH,\m)\to (\cH',\m')$ 
and the notion of $A_\infty$-quasi-isomorphisms in fact defines 
a homotopy equivalence relation between $A_\infty$-algebras 
(see \cite{hk:thesis} and reference therein).

 \subsection{Homological perturbation theory for $A_\infty$-structures}
\label{ssec:HPT}

A version of homological perturbation theory we shall employ 
is as follows. 
\begin{thm}\label{thm:HPT}
For an $A_\infty$-algebra $(\cH,\m)$, suppose given linear maps 
$h:\cH^r\to\cH^{r-1}$ and $P:\cH^r\to\cH^r$ 
satisfying 
\begin{equation}\label{HKdecomp}
 dh+hd=\Id_{\cH}-P, \quad P^2=P, \qquad d:=m_1
\end{equation}
on $\cH$. Then, there exists a canonical way to construct 
an $A_\infty$-structure $\m'$ on $P\cH$ such that 
$(P\cH,\m')$ is homotopy equivalent 
to the original $A_\infty$-algebra $(\cH,\m)$. 
\end{thm}
Note that if $d P=0$, then eq.(\ref{HKdecomp}) gives a Hodge 
decomposition of the complex $(\cH,d)$, where 
$P(\cH)=H(\cH)$ gives the cohomology. 
\begin{pf}
Let $\iota:\cH'\to \cH$ be the embedding and $\pi:\cH\to\cH'$ 
the projection, respectively, such that 
$\pi\circ\iota=\Id_{\cH'}$ and $\iota\circ\pi=P$. 
Namely, $\iota$ is the embedding $\cH'\simeq P(\cH)\subset\cH$. 
A collection of degree zero maps 
$\cG=\{g_l:(\cH')^{\otimes l}\to\cH\}_{l\ge 1}$ 
is defined recursively 
with respect to $k$ as 
 \begin{equation}\label{morp-recursive}
 g_k =-h\sum_{i\ge 2}\sum_{1\le k_1<k_2\cdots<k_i=k}
 m_i(g_{k_1}\otimes g_{k_2-k_1}\otimes
 \cdots\otimes g_{k-k_{i-1}})
 \end{equation}
with $g_1:=\iota:\cH'\to\cH$ the inclusion. 
Then, $\m'=\{m'_k:(\cH')^{\otimes k}\to\cH\}_{k\ge 1}$ 
is given recursively by 
 \begin{equation}\label{prod-recursive}
  m'_k = \pi\sum_{i\ge 2}\sum_{1\le k_1<k_2\cdots <k_i=k}
 m'_i(g_{k_1}\otimes g_{k_2-k_1}\otimes
 \cdots\otimes g_{k-k_{i-1}})\ .
 \end{equation}
Note that $(m_1)^2=\pi\circ d\circ\iota\circ\pi\circ d\circ\iota
=\pi\circ d\circ P\circ d\circ\iota=0$ 
since $d$ commutes with $P$ due to the condition (\ref{HKdecomp}). 
One can check that 
these actually give an $A_\infty$-structure and 
an $A_\infty$-quasi-isomorphism 
(see \cite{hk:thesis}). 
\qed\end{pf}
Equivalently, $\m'$ are described in terms of rooted planar trees 
as follows. 

A {\em planar tree} (a simply connected 
planar graph without loops) 
consists of vertices, internal edges and external edges. 
An internal edge has two distinct vertices at its ends. 
An external edge has one end on a vertex and another end is free. 
The number of incident edges at a vertex is greater than two. 
The term `planar' means the cyclic order of edges at each vertex is 
distinguished. 
A {\em rooted planar tree} is a planar tree graph with one of its external 
edges distinguished from others as a root edge. 
The remaining external edges are called the {\em leaves}. 
Each edge of a planar rooted tree has a unique orientation 
so that the orientations form a flow from the leaves to the root edge. 
We sometimes describe the orientation as an arrow. 
We call a vertex at which the number of incident edges is $(k+1)$ 
a {\em $k$-vertex}.

We call a rooted planar tree having $k$ leaves a {\em $k$-tree}. 
The set of (the isomorphism classes of) $k$-trees is 
denoted by $G_k$, $k\ge 2$. 

For any element $\Gamma_n\in G_n$, $n\ge 2$, 
let us define $m'_{\Gamma_n}:(\cH')^{\otimes n}\to\cH'$ by 
attaching $\iota:\cH'\to\cH$ to each leaf, $m_k:\cH^{\otimes k}\to\cH$ to 
each $k$-vertex, $-h:\cH\to\cH$ to 
each internal edge, $\pi:\cH\to\cH'$ to the root edge and 
then composing them. 
For example, 
\begin{equation*}
 \begin{minipage}[c]{70mm}{
 \begin{equation*}
 \begin{split}
 & m'_{\Gamma_3}(o'_1,o'_2,o'_3) \\
 & \ \ \ = \pi m_2(-hm_2(\iota(o'_1),\iota(o'_2)),\iota(o'_3))
 \end{split}
 \end{equation*}
 }\end{minipage}, \qquad \Gamma_3=
 \begin{minipage}[c]{50mm}{\scalebox{0.5}[0.5]{\includegraphics{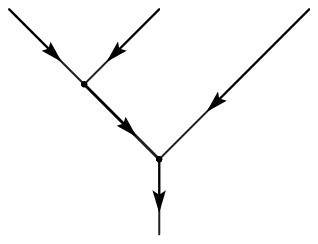}}}\end{minipage}
\end{equation*}
for $o'_1,o'_2,o'_3\in\cH'$. 
Then, $\{m'_n\}_{n\ge 1}$ is given by $m'_1=\pi\circ m_1\circ\iota$ and 
\begin{equation}\label{HPTtree}
 m'_n=\sum_{\Gamma_n\in G_n} m'_{\Gamma_n}
\end{equation}
for $n\ge 2$. Thus, $m'_n$ is described as the sum of the value 
$m'_{\Gamma_n}$ over all the $n$-trees $\Gamma_n\in G_n$. 
Similarly, $\{g_n\}_{n\ge 1}$ is given by $g_1=\iota$ and 
$g_n=\sum_{\Gamma_n\in G_n}g_{\Gamma_n}$ for $n\ge 2$, 
where $g_{\Gamma_n}:(\cH')^{\otimes n}\to\cH$ is obtained by 
replacing $\pi$ by $-h$ in the definition of $m'_{\Gamma_n}$. 
\begin{rem}\label{rem:HPT}
The data $\Nsddata{P\cH}{\ \iota}{\ \pi}{\cH}{h}$ 
used in the proof above is often called a 
{\em strong deformation retract (SDR)} of the complex $(\cH,d)$, 
the starting point of the traditional HPT (for instance 
\cite{gugjds,GLS:chen,GLS,HK}). 
There, it is discussed that 
the $A_\infty$-quasi-isomorphism (\ref{morp-recursive}) induces 
homotopy equivalence of 
the induced $A_\infty$-products (\ref{prod-recursive}) 
with the original one $\m$, 
mainly in the case $(\cH,\m)$ is a DGA. 
The extension to the case when $(\cH,\m)$ is a general $A_\infty$-algebra 
is not difficult. 
The present form of HPT (Theorem \ref{thm:HPT}) is due to \cite{KoSo}, 
where the above planar tree expression of the recursive formula 
eq.(\ref{morp-recursive}) and (\ref{prod-recursive}) is also presented. 
\end{rem}

 \subsection{$A_\infty$-categories}
\label{ssec:Ainfty-cat}

We need the categorical version of these terminologies. 
\begin{defn}[$A_\infty$-category \cite{fukayaAinfty}]
An {\em $A_\infty$-category} $\cC$ consists of 
the set of objects $\Ob(\cC)=\{a,b,\dots\}$, 
$\Z$-graded vector space $V_{ab}:=\Hom_\cC(a,b)$ for 
each two objects $a,b\in\Ob(\cC)$ and a collection of 
multilinear maps 
\begin{equation*}
 \m:=\{m_n:V_{a_1a_2}\otimes\cdots\otimes V_{a_na_{n+1}}
 \to V_{a_1a_{n+1}}\}_{n\ge 1}
\end{equation*}
of degree $(2-n)$ defining the $A_\infty$-structure, that is, 
$\m$ satisfies the $A_\infty$-relations 
(\ref{Ainfty}). 

In particular, an $A_\infty$-category $\cC$ with 
vanishing higher products $m_3=m_4=\cdots=0$ is called 
a {\em DG category}. 
\end{defn}
The suspension $s(\cC)$ of an $A_\infty$-category $\cC$ is defined by 
the shift 
\begin{equation*}
 s: \Hom_{\cC}(a,b)\to s(\Hom_{\cC}(a,b))=:\Hom_{s(\cC)}(a,b)
\end{equation*}
for any $a,b\in\Ob(\cC)=\Ob(s(\cC))$, 
where the degree $|m_n|$ of the $A_\infty$-products becomes one 
for all $n\ge 1$ as in the case of $A_\infty$-algebras. 
We sometimes denote $\Hom_{s(\cC)}(a,b)=\cH_{ab}$ as 
we do $\Hom_{\cC}(a,b)=V_{ab}$.

\begin{defn}[$A_\infty$-functor]
Given two $A_\infty$-categories $\cC$, $\cC'$, 
$\cG:=\{g,g_1,g_2,\dots\}:s(\cC)\to s(\cC')$ is called an 
{\em $A_\infty$-functor} if and only if 
$g:\Ob(s(\cC))\to\Ob(s(\cC'))$ is a map of object and 
\begin{equation*}
 g_k:\Hom_{s(\cC)}(a_1,a_2)\otimes\cdots\otimes\Hom_{s(\cC)}(a_k,a_{k+1})
 \to\Hom_{s(\cC')}(g(a_1),g(a_{k+1})),\quad k\ge 1
\end{equation*}
are degree preserving multilinear maps satisfying the defining relations 
of an $A_\infty$-morphism (\ref{Ainftymorp}). 

In particular, if $g:\Ob(s(\cC))\to\Ob(s(\cC'))$ and 
$g_1:\Hom_{s(\cC)}(a,b)\to\Hom_{s(\cC')}(f(a),f(b))$ 
induces an isomorphism between the 
cohomologies for any $a,b\in\Ob(s(\cC))$, 
we call the $A_\infty$-functor {\em homotopy equivalence}. 
\end{defn}
The generalization of HPT for $A_\infty$-algebras to 
$A_\infty$-categories is straightforward \cite{KoSo}. 
\begin{thm}\label{thm:HPT-cat}
For an $A_\infty$-category $\cC$, 
suppose given linear maps 
$h_{ab}:\cH^r_{ab}\to\cH^{r-1}_{ab}$ and $P_{ab}:\cH^r_{ab}\to\cH^r_{ab}$ 
satisfying 
\begin{equation}\label{HKdecomp2}
 d_{ab}h_{ab}+h_{ab}d_{ab}=\Id_{\cH_{ab}}-P_{ab}, \quad (P_{ab})^2=P_{ab}, 
 \qquad d_{ab}:=m_1:\cH_{ab}\to\cH_{ab}
\end{equation}
on $\cH_{ab}$ for any $a,b\in\Ob(\cC)$. 
Then, there exists a canonical way to construct 
an $A_\infty$-category $\cC'$ which is homotopy equivalent to the original 
$A_\infty$-category $\cC$ and in particular 
the space of morphisms is defined by $\Hom_{s(\cC')}(a,b)=\cH'_{ab}=P_{ab}\cH_{ab}$. 
\end{thm}
\begin{pf}
Let $\iota_{ab}:\cH'_{ab}\to\cH_{ab}$ be the embedding 
and $\pi_{ab}:\cH_{ab}\to\cH'_{ab}$ 
the projection such that $\pi_{ab}\circ\iota_{ab}=\Id_{\cH'_{ab}}$ and 
$\iota_{ab}\circ\pi_{ab}=P_{ab}$. Then, 
for $a_1,\dots,a_{n+1}\in\Ob(\cC')$, the $A_\infty$-product 
$m_n':\cH'_{a_1a_2}\otimes\cdots\otimes\cH'_{a_na_{n+1}}\to\cH'_{a_1a_{n+1}}$ 
is given by $m'_n=\sum_{\Gamma_n\in G_n}m'_{\Gamma_n}$, 
where $m'_{\Gamma_n}$ is defined in the same way as the one 
for an $A_\infty$-algebra, 
but we attach $\iota_{a_ia_{i+1}}:\cH'_{a_ia_{i+1}}\to\cH_{a_ia_{i+1}}$, $i=1,\dots,n$, 
for each leaf (instead of $\iota$), 
$m_k$ to each $k$-vertex, 
$h_{ab}$ to each internal edge, where $a,b\in\{a_1,\dots,a_{n+1}\}$ 
is uniquely determined by the graph $\Gamma_n$, and finally 
$\pi_{a_1a_{n+1}}$ to the root edge of $\Gamma_n$ (instead of 
$\pi$). The construction of homotopy equivalence is also parallel to the case of 
$A_\infty$-algebras, though we do not use it in the present paper. 
\qed\end{pf}

 \section{$A_\infty$-category of lines in a plane}

 \subsection{The main theorem}
\label{ssec:main}

For a fixed integer $N\ge 2$, 
let $\{f_1,\dots, f_N\}$ be a set of polynomial 
functions on $\R$ of degree equal or less than two. 
For each $a\in\{1,\dots,N\}$, 
$y= df_a/dx$ is a line $L_a$ in $\R^2$ with coordinates 
$(x,y)$ described as 
\begin{equation*}
 L_a : y=t_a x + s_a,\qquad  t_a,s_a\in\R. 
\end{equation*} 
Let us consider such a collection $\{f_1,\dots,f_N\}$ 
satisfying the following two conditions: 
\begin{itemize}
 \item[(i)] For any $a\ne b=1,\dots,N$, the slopes of the 
lines $L_a$ and $L_b$ are not the same: $t_a\ne t_b$. 

 \item[(ii)] More than two lines do not intersect at the same point in $\R^2$. 
\end{itemize}
We identify the set $\{f_a\ |a=1,\dots,N\}$ with the label set 
$\{a | a=1,\dots,N\}$. 
Then, denote by $\fF_N:=\{a | a=1,\dots,N\}$ such a set 
satisfying the above conditions (i) and (ii). 

We shall construct a Fukaya $A_\infty$-category $\cC(\fF_N)$ with 
$\Ob(\cC(\fF_N))=\fF_N$ from another $A_\infty$-category, in particular, 
a DG category $\cC_{DR}(\fF_N)$. 
Let $\Omega(\R):=\Omega^0(\R)\oplus\Omega^1(\R)$ be 
the graded vector space defined by 
$\Omega^0(\R):=\cS(\R)$, the space of Schwartz functions, and 
$\Omega^1(\R):=\cS(\R)\cdot dx$, where $dx$ is the base of one-form on $\R$. 
\begin{defn}[$\cC_{DR}(\fF_N)$]
The DG category $\cC_{DR}(\fF_N)$ consists of the set of objects 
$\Ob(\cC_{DR}(\R))=\fF_N$ and 
the space of morphisms $\Omega_{ab}:=\Hom_{\cC_{DR}(\fF_N)}(a,b)=\Omega(\R)$ 
for each $a,b\in\fF_N$, 
where we set 
\begin{itemize}
 \item the differential $d_{ab}:\Omega_{ab}^0\to\Omega^1_{ab}$ 
by $d_{ab}:=d-df_{ab}\wedge$, where $d=dx\cdot (d/dx)$ 
is the exterior derivative and 
$f_{ab}:=f_a-f_b$; 

 \item the product 
$m: \Omega^{r_{ab}}_{ab}\otimes\Omega^{r_{bc}}_{bc}
\to\Omega^{r_{ab}+r_{bc}}_{ac}$ 
by the usual wedge product. 
\end{itemize}
\end{defn}
It is clear that $\cC_{DR}(\fF)$ forms a DG category. 

The following is the main theorem of this paper. 
\begin{thm}\label{thm:main}
There exists an $A_\infty$-category $\cC(\fF_N)$ with 
$\Ob(\fF_N)=\fF_N$ such that 
\begin{itemize}
 \item[(i)] For two objects $a\ne b\in \fF_N$, the space 
$\Hom_{\cC(\fF_N)}(a,b)=:V_{ab}$ of morphisms is 
the following graded vector space of degrees zero and one$:$ 
\begin{equation*}
\renewcommand{\arraystretch}{1.3}
 \begin{array}{lll}
   V^0_{ab}=\R\cdot [v_{ab}], &V^1_{ab}=0,  &\quad t_a< t_b,  \\
   V^0_{ab}=0, & V^1_{ab}=\R\cdot [v_{ab}], &\quad t_a> t_b. 
 \end{array}
\end{equation*}
Here, $[v_{ab}]$ are the bases 
of the vector spaces attached to the intersection points 
$v_{ab}(=v_{ba})$ of $L_a$ and $L_b$. 

 \item[(ii)]
Let $a_1,\dots,a_{k+1}\in\fF_N$, $k\ge 1$, be objects such that $a_i\ne a_j$ 
for any $i\ne j\in\{1,\dots,k+1\}$ and 
$\vec{v}:=(v_{a_1},\dots,v_{a_ka_{k+1}},v_{a_{k+1}a_1})$. 
Then, for $k=1$, the differential $m_1:V_{a_1a_2}\to V_{a_1a_2}$ 
is zero, $m_1=0$. For $k\ge 2$, 
the structure constant $c(\vec{v})\in\R$ for 
the higher $A_\infty$-product 
\begin{equation*}
 m_k([v_{a_1a_2}],\dots,[v_{a_ka_{k+1}}])
 = c(\vec{v})\cdot [v_{a_1a_{k+1}}] 
\end{equation*}
is zero if $\vec{v}$ does not form a clockwise convex polygon 
(see also Definition \ref{defn:polygon} 
for the definition of clockwise convex polygon), 
and if $\vec{v}$ forms a clockwise convex polygon, it is 
given by $c(\vec{v})=\pm e^{-Area(\vec{v})}$, 
with an appropriate sign $\pm$, where $Area(\vec{v})$ is the area of 
the clockwise convex polygon. 

 \item[(iii)] $\cC(\fF_N)$ is homotopy equivalent to $\cC_{DR}(\fF_N)$. 
 \end{itemize}
\end{thm}
Conditions (i) and (ii) are the ones for $\cC(\fF_N)$ to be a 
Fukaya category. 
We call a multilinear map $m_k$, $k\ge 2$, of the type in Condition (ii) 
a {\em transversal (higher) $A_\infty$-product}. 
Multilinear maps $m_k$ of the other type are then called 
{\em non-transversal $A_\infty$-products}. 
Condition (iii) is motivated by homological mirror symmetry (HMS)\cite{HMS} 
of (non)commutative complex tori. 
As discussed in \cite{KoSo}, 
this homotopy equivalence should be the key idea of HMS 
for tori or more general cases, where 
both $\cC(\fF_N)$ and $\cC_{DR}(\fF_N)$ are $A_\infty$-categories 
associated to a symplectic structure, but 
$\cC_{DR}(\fF_N)$ is canonically isomorphic to a DG category 
associated to the mirror dual complex structure. 
In fact, the relation of this DG category $\cC_{DR}(\fF_N)$ 
with the DG category of holomorphic vector bundles on 
a noncommutative complex torus \cite{PoSc, hk:nchms} with 
noncommutativity set to be zero is clear. 
For the precise relation of the {\em noncommutative} 
complex torus description 
and the usual complex torus description, for instance see \cite{hk:nctheta}.

 \subsection{Commutative DG algebra $A_S(F)$}
\label{ssec:CDGA}

We shall give a sketch of the proof of Theorem \ref{thm:main} 
in section \ref{sec:proof}. 
Before that, we present an $A_\infty$-category, 
which hereafter we denote by $\cC(\fF_N)$, 
shown to satisfy Conditions (i), (ii) and (iii) in the next subsection. 

In order to construct an $A_\infty$-structure 
including non-transversal $A_\infty$-products, 
we introduce a (commutative) DG algebra $A_S(F)$ over $k=\R$. 
This notion is motivated by 
an extension of a subalgebra $F$ of the commutative DG algebra 
of smooth differential forms on $\R$ by including 
step functions and delta-function one forms. 
\begin{defn}[Commutative DG algebra $A_S(F)$]
Let $F=F^0\oplus F^1$ be a commutative DG subalgebra 
of the commutative DG algebra of smooth differential forms on $\R$, 
and $S$ be a finite set with a map $x: S\to\R$. 
For each $v\in S$, we introduce degree zero base $\vartheta_v$ 
and degree one base $\delta_v:=d(\vartheta_v)$ 
with (degree zero) unit $\1$: 
$\1\cdot\vartheta_v=\vartheta_v\cdot\1=\vartheta_v$, 
$\1\cdot\delta_v=\delta_v\cdot \1=\delta_v$. 
Consider the commutative algebra 
$\ti{A}_S(F):=F\otimes\la \1,\vartheta_v,\delta_v\ |\ v\in S\ra$ 
of degrees zero and one, and relations defined as follows: 
\begin{equation*} 
 \vartheta_v\vartheta_{v'}=\vartheta_{v'},\quad \delta_v\vartheta_{v'}=0
\end{equation*}
for any $v,v'\in S$ such that $x(v)< x(v')$, 
\begin{equation*} 
 \vartheta_v=\vartheta_{v'},\quad \delta_v=\delta_{v'}
\end{equation*}
for any $v,v'\in S$ such that $x(v)= x(v')$, 
\begin{equation*}
 \alpha\cdot\delta_v=\alpha(x(v))\cdot\delta_v,\qquad \alpha(x(v))\in k=\R, 
\end{equation*}
for any $v\in S$, 
\begin{equation*}
 \alpha\cdot\vartheta_v=0,\quad 
 \beta\cdot\vartheta_v=0,\qquad \alpha\in F^0,\ \ \beta\in F^1
\end{equation*}
for any $v\in S$ if $\alpha(x)=0$ or $\beta(x)$=0 for any $x\ge x(v)$, 
\begin{equation*}
 \alpha\cdot (\1-\vartheta_v)=0,\quad 
 \beta\cdot (\1-\vartheta_v)=0,\qquad \alpha\in F^0,\ \ \beta\in F^1
\end{equation*}
for any $v\in S$ if $\alpha(x)=0$ or $\beta(x)=0$ for any $x\le x(v)$, 
$F^1\cdot\delta_v=0$ for any $v\in S$ and 
$\delta_v\cdot\delta_{v'}=0$ for any $v,v'\in S$. 
More explicitly, 
the graded vector space $\ti{A}^r_S(F)$, $r=0,1$, is 
\begin{equation*}
 \ti{A}_S^0(F):=F^0\otimes\la \1,\vartheta_v | v\in S\ra,\qquad 
 \ti{A}_S^1(F):=F^1\otimes \la \1,\vartheta_v | v\in S\ra
 \oplus \oplus_{v\in S}\ti{A}^0_S(F)\otimes\delta_v. 
\end{equation*}
By the commutativity and the relations above, 
any element is described as 
\begin{equation*}
 \alpha=\alpha_0 + \sum_{v\in S, n\in\Z_{>0}}
 \alpha_{v,n}(\vartheta_v)^n,\qquad 
 \alpha_0,\ \alpha_{v,n}\in F^0
\end{equation*}
for $\alpha\in\ti{A}_S^0(F)$ and 
\begin{equation*}
 \beta=\beta_0 + \sum_{v\in S, n\in\Z_{>0}}\beta_{v,n}(\vartheta_v)^n
 + \sum_{v\in S, n\in\Z_{> 0}} c_{v,n}(\vartheta_v)^{n-1}\cdot\delta_v,\qquad 
 \beta_0,\ \beta_{v,n}\in F^1,\quad c_{v,n}\in k=\R 
\end{equation*}
for $\beta\in\ti{A}_S^1(F)$. 
The differential $d: \ti{A}_S^0(F)\to\ti{A}_S^1(F)$ is 
defined by extending the differential $d:F^0\to F^1$ 
with $d(\vartheta_v)=\delta_v$, $v\in S$, so that 
they satisfy the Leibniz rule with respect to the commutative product. 

In this paper, we shall consider 
the two cases $F=\Omega(\R)$ and $F=F^0=\R$. 
For $F=\Omega(\R)$, we set $A_S(\Omega(\R)):=\ti{A}_S(\Omega(\R))$. 
For $F=\R$ (note that the differential on $F$ is trivial), 
we set $A_S(\R)$ as a commutative DG subalgebra of $\ti{A}_S(F)$ as follows: 
\begin{equation*}
 \begin{split}
 & A^0_S(\R):=
\left\{\alpha_0+\sum_{v\in S, n\in\Z_{>0}}\alpha_{v,n}
 (\vartheta_v)^n\in\ti{A}_S(\R) 
 \Big|\ \alpha_0=0, \sum_{v\in S, n\in\Z_{>0}}\alpha_{v,n}=0 \right\}\ ,\\
 & A^1_S(\R):=\ti{A}^1_S(\R)
=\left\{ \sum_{v\in S, n\in\Z_{> 0}} c_{v,n}
 (\vartheta_v)^{n-1}\cdot\delta_v,\ \ 
 c_{v,n}\in k=\R\ \right\}. 
 \end{split}
\end{equation*}
\end{defn}
By the map $x:S\to\R$, $S$ is identified 
with the set of finitely many points on $\R$. Then, 
$\delta_v$ can be regarded as the delta function on $\R$ 
with support at $x(v)$ and $\vartheta_v$ is the step function whose value is zero at 
$x\to -\infty$, one at $x\to\infty$ and which is discontinuous 
at $x(v)$. In the case $F=\R$, 
taking the subspace $A^0_S(\R)\subset\ti{A}^0_S(\R)$ implies that 
we concentrate on constant functions $\alpha$ which are discontinuous 
at some points in $x(S)$ and further $\alpha=0$ at $x\to\pm\infty$. 
Note that we do not impose the relation $(\vartheta_v)^2=\vartheta_v$. 
Thus, for any element in $A^0_S(\R)$, any discontinuous point $x(v)$ 
is associated with $\Z_{> 0}$-valued weight 
corresponding to the power of $\vartheta_{v}$. 
The commutative DGA 
$A_S(\R)$ is useful in the sense that it is defined only 
in terms of finitely many points on $\R$, though 
$A_S(\R)$ is infinite dimensional as a vector space. 

For the construction of the $A_\infty$-category $\cC(\fF_N)$ 
we need only $A_S(\R)$. 
The cohomology of $A_S(\R)$ is 
$H^0(A_S(\R))=0$ and $H^1(A_S(\R))\simeq\R$ (one dimensional); 
a base of $H^1$ is $\delta_v$ for an element $v\in S$, 
but one has $\delta_{v'}-\delta_v=d(\vartheta_{v'}-\vartheta_v)$ 
for $v,v'\in S$.

At a first look 
the reader can skip the following lemma, 
which shall be employed as a key step of the proof of 
Theorem \ref{thm:main} in subsection \ref{ssec:proof}. 
\begin{lem}\label{lem:CDG}
There exist inclusions 
\begin{equation*}
 \iota: A_S(\R) \to A_S(\Omega(\R)),\qquad 
 \iota: \Omega(\R) \to A_S(\Omega(\R)),
\end{equation*}
both of which induce homotopy equivalences as $A_\infty$-algebras. 
\end{lem}
\begin{pf}
The existence of the inclusions $\iota$ is clear. Also, for each case, 
the $\iota$ forms a chain map with respect to the differentials 
on both sides, 
and also defines an algebra homomorphism. 
Then, for each case, by setting $g_1:=\iota$ and $g_2=g_3=\cdots=0$, 
$\cG:=\{g_1,g_2,\dots\}$ forms an $A_\infty$-morphism. 

For $A_S(\R)$, $\Omega(\R)$ and $A_S(\Omega(\R))$, their cohomologies are 
isomorphic to each other: $H^0=0$ and $H^1=\R$. 
In order to show that $g_1:=\iota$ induces isomorphism on the cohomologies, 
we need only see that the image of a representative of $H^1$ of 
$A_S(\R)$ or $\Omega(\R)$ is not exact in $A_S(\Omega(\R))$. 
It is clear that the image of $\delta_v\in A_S(\R)$ and 
the image of $\beta_0\in\Omega^1(\R)$ such that 
$\int_{-\infty}^{\infty}\beta_0\ne 0$ are not exact.
\qed\end{pf}

 \subsection{The $A_\infty$-category $\cC(\fF_N)$}
\label{ssec:explicit}

Let us define the $A_\infty$-category $\cC(\fF_N)$. 
First of all, for each $a\in\fF_N$ 
the graded vector space $V_{aa}=V_{aa}^0\oplus V_{aa}^1$ 
is set to be 
\begin{equation*}
 V^r_{aa}=A_{S_a}^r(\R),\qquad r=0,1, 
\end{equation*}
on which we set the differential $m_1=d:V_{aa}^0\to V_{aa}^1$ 
and the product $m_2:V_{aa}\otimes V_{aa}\to V_{aa}$ 
as those in $A_{S_a}(\R)$. 

For $a\ne b\in\fF_N$, 
the graded vector space $V_{ab}$ is taken to be the one given 
in Theorem \ref{thm:main} (i), on which the differential 
$m_1:V_{ab}\to V_{ab}$ is set to be zero. 
Then, next, let us define multilinear maps 
\begin{equation*}
 m_k: V_{a_1a_2}\otimes\cdots\otimes V_{a_ka_{k+1}}\to V_{a_1a_{k+1}} 
\end{equation*}
of degree $(2-k)$ for $k\ge 2$. By degree counting, the following holds. 
\begin{lem}\label{lem:degree}
Any multilinear maps $m_k(w_1,\dots,w_k)$ can be nonzero only if 
there exists a nonzero element $w_{k+1}\in V_{a_{k+1}a_1}$ 
such that the number of degree zero elements 
in $\{w_1,\dots,w_{k+1}\}$ is two. 
\qed\end{lem}
We first define multilinear maps $m_k(w_1,\dots,w_k)$ on $\ti{V}_{**}$, 
where $\ti{V}_{ab}=V_{ab}$ for $a\ne b\in\fF_N$ 
and $\ti{V}_{aa}=\ti{A}_{S_a}(\R)$ for $a\in\fF_N$. 
The multilinear maps on the $\Z$-graded vector spaces $\ti{V}_{**}$ given below 
are closed in the $\Z$-graded subvector spaces $V_{**}$ 
and thus the restriction of them onto $V_{**}$ 
gives the multilinear maps on $V_{**}$. 

We determine those multilinear maps on $\ti{V}_{ab}$ 
separately in each case 
$\sharp(\vartheta):=\sharp\{1\le i\le k | w_i\in V_{aa}^0,\ a\in\fF_N\}$ 
is two, one, or zero. 

\vspace*{0.3cm}

\noindent
$\bullet$\ {\bf The case $\sharp(\vartheta)=2$}: \ \ 
By degree counting (Lemma \ref{lem:degree}), 
the multilinear map $m_k(w_1,\dots,w_k)$ can be nonzero only if 
$w_i\in\ti{V}_{aa}$ for all $i=1,\dots, k$ with some $a\in\fF_N$. 
We set $m_k(w_1,\dots,w_k)$ is nonzero only if it is of the form 
$m_2(w_1,w_2)$, $w_1,w_2\in\ti{V}_{aa}^0$, for some $a\in\fF_N$. 
This is the product $m_2$ in $\ti{A}_{S_a}(\R)$. 

\vspace*{0.3cm}

\noindent
$\bullet$\ 
{\bf The case $\sharp(\vartheta)=1$}: \ \ 
By degree counting (Lemma \ref{lem:degree}), 
all such $A_\infty$-products to be nonzero are only of the following types: 
for any $a\ne b\in\fF_N$, 
\begin{itemize}
 \item[A$_{\ }$] \quad 
 $
 m_*: (V_{aa}^1)^{\otimes k_1}\otimes 
 \ti{V}_{aa}^0\otimes (V_{aa}^1)^{\otimes k_2}
 \otimes V_{ab}^r\otimes (V_{bb}^1)^{\otimes k_3}\to V_{ab}^r$, 
 \item[B$_{\ }$] \quad 
 $
 m_*: (V_{aa}^1)^{\otimes k_1}\otimes V_{ab}^r
 \otimes (V_{bb}^1)^{\otimes k_2}\otimes \ti{V}_{bb}^0
 \otimes (V_{bb}^1)^{\otimes k_3}\to V_{ab}^r$, 
 \item[C$_{1}$]\quad 
 $
m_*: 
 (V_{aa}^1)^{\otimes k_1}\otimes V_{ab}^r\otimes (V_{bb}^1)^{\otimes k_2}
 \otimes\ti{V}_{bb}^0\otimes (V_{bb}^1)^{\otimes k_3}
 \otimes V_{ba}^{1-r}\otimes (V_{aa}^1)^{\otimes k_4}\to\ti{V}_{aa}^0$, 

 \item[C$_{2}$]\quad 
 $
 m_*: (V_{aa}^1)^{\otimes k_1}\otimes\ti{V}_{aa}^0 
 \otimes (V_{aa}^1)^{\otimes k_2}
 \otimes V_{ab}^1\otimes (V_{bb}^1)^{\otimes k_3}
 \otimes V_{ba}^{0}\otimes (V_{aa}^1)^{\otimes k_4}\to\ti{V}_{aa}^0$, 

 \item[C$_{3}$]\quad 
 $
 m_*: (V_{aa}^1)^{\otimes k_1}
 \otimes V_{ab}^0\otimes (V_{bb}^1)^{\otimes k_2}\otimes V_{ba}^1
 \otimes (V_{aa}^1)^{\otimes k_3}
 \otimes\ti{V}_{aa}^{0}\otimes (V_{aa}^1)^{\otimes k_4}\to\ti{V}_{aa}^0$, 
\end{itemize}
where $*$'s are the appropriate numbers. 
We set the $A_\infty$-products of the following types to be zero; 
type A with $r=0$ if $k_1\ne 0$, 
type A with $r=1$ if $k_2\ne 0$, 
type B with $r=0$ if $k_3\ne 0$, 
type B with $r=1$ if $k_2\ne 0$, 
type C$_{1}$ with $r=0$ if $k_3\ne 0$, 
type C$_{1}$ with $r=1$ if $k_2\ne 0$, 
type C$_{2}$ and type C$_{3}$ if $k_2\ne 0$. 

We set the multilinear maps which do not include 
degree one elements in $V_{aa}^1$ for any $a\in\fF_N$ as follows. 

For type A, the product $m_2:\ti{V}_{aa}^0\otimes V_{ab}^r\to V_{ab}^r$, 
$a\ne b$, $r=0, 1$, is given by 
\begin{equation}\label{l-module}
 m_2((\vartheta_{v_a})^n,[v_{ab}])
 \begin{cases}
 [v_{ab}] & \quad x(v_a) < x(v_{ab}) \\
 \ov{2^n}[v_{ab}] & \quad v_a=v_{ab},\quad t_a<t_b \\
 \ov{n+1}[v_{ab}] & \quad v_a=v_{ab},\quad t_a>t_b \\
 0 & \quad  x(v_{ab})< x(v_a) 
 \end{cases}
\end{equation}
for $n\ge 1$, where recall that the degree of 
$[v_{ab}]$ is zero for $t_a<t_b$ and one for $t_a>t_b$. 
In the same way, for type B, 
the product $V_{ab}^r\otimes\ti{V}_{bb}^0\to V_{ab}^r$, 
$a\ne b$, $r=0,1$, is given by 
\begin{equation}\label{r-module}
 m_2([v_{ab}],(\vartheta_{v_b})^n)=
 \begin{cases}
 [v_{ab}] & \quad x(v_b) < x(v_{ab}) \\
 \ov{2^n}[v_{ab}] & \quad v_b=v_{ab},\quad t_a<t_b \\
 \ov{n+1}[v_{ab}] & \quad v_b=v_{ab},\quad t_a>t_b \\
 0 & \quad  x(v_{ab})< x(v_b) 
 \end{cases}
\end{equation}
for $n\ge 1$. 
In addition, we set $m_2(\1_a,[v_{ab}])=[v_{ab}]$ 
and $m_2([v_{ab}],\1_b)=[v_{ab}]$ 
for the identities $\1_a\in\ti{V}_{aa}^0$ and $\1_b\in\ti{V}_{bb}^0$.

For type C$_{1}$, C$_{2}$, C$_{3}$, $a\ne b\in\fF_N$ such that $t_a<t_b$, 
\begin{equation*}
 \begin{split}
 & m_3([v_{ba}],(\vartheta_{v_a})^n,[v_{ab}]) =
 \ov{n+1}\vartheta_{v_{ab}}(1-(\vartheta_{v_{ab}})^n)\ \in\ti{V}_{bb}^0, \\
 & m_3((\vartheta_{v_b})^n,[v_{ba}],[v_{ab}]) = 
 \ov{n+1}\vartheta_{v_{ab}}(1-(\vartheta_{v_{ab}})^n)\ \in\ti{V}_{bb}^0, \\
 & m_3([v_{ab}],(\vartheta_{v_b})^n,[v_{ba}]) = 
 -\ov{n+1}\vartheta_{v_{ab}}(1-(\vartheta_{v_{ab}})^n)\ \in\ti{V}_{aa}^0, \\
 & m_3([v_{ab}],[v_{ba}],(\vartheta_{v_a})^n) = 
 -\ov{n+1}\vartheta_{v_{ab}}(1-(\vartheta_{v_{ab}})^n)\ \in\ti{V}_{aa}^0 
 \end{split}
\end{equation*}
for $n\ge 1$ if $v_a=v_{ab}$ or $v_b=v_{ab}$, 
and they are equal to zero if $v_a\ne v_{ab}$ or $v_b\ne v_{ab}$. 
In addition, we set 
$m_3([v_{ba}],\1_a,[v_{ab}])=m_3(\1_b,[v_{ba}],[v_{ab}])=0$ 
and 
$m_3([v_{ab}],\1_b,[v_{ba}])=m_3([v_{ab}],[v_{ba}],\1_a)=0$.

\vspace*{0.3cm}

\noindent
$\bullet$\ {\bf The case $\sharp(\vartheta)=0$}:\ \ 
We first prepare some terminology for polygons. 
\begin{defn}[CC-polygon, Degree of points, Sign of the CC-polygon]
\label{defn:polygon}
Let $\vec{v}$ be a sequence of points $v_{ab}$, 
$a,b\in\fF_N$, in $\R^2$ with coordinates $(x,y)$. 
Any $\vec{v}$ is described in the form 
\begin{equation*}
 \vec{v}=(v_1,\dots,v_1,v_2,\dots,v_2,\dots\dots, v_n,\dots,v_n), 
\end{equation*}
where $\{v_1,\dots,v_n\}$, $n\in\Z_{> 0}$, are points in $\R^2$ 
such that $v_i\ne v_{i+1}$ for $i=1,\dots, n-1$. 
In this expression, 
we call $\vec{v}$ a {\em point} if $n=1$. 
On the other hand, we call $\vec{v}$ a 
{\em clockwise convex polygon (CC-polygon)} 
if and only if $0< Angle(v_{i-1}v_iv_{i+1})\le \pi$ for any $i\in\Z$, 
where, we identify $v_i=v_{i+n}$ 
if $v_1\ne v_n$ and $v_i=v_{i+(n-1)}$ if $v_1 = v_n$.  
By definition $n\ge 3$ if $\vec{v}$ is a CC-polygon. 

For a CC-polygon 
$\vec{v}=(v_1,\dots,v_1,v_2,\dots,v_2,\dots\dots, v_n,\dots,v_n)$, 
we attach a {\em degree} $|v_i|$ for each point $v_i$, $i=1,\dots,n$, 
as follows. 
Consider the map $x:\{v_1,\dots, v_n\}\to\R$, where 
the image $x(v_i)$ is the $x$-coordinate of the point $v_i$. 
Let $\{x_L<\cdots< x_R\}\subset\R$ be the ordered subset 
consisting of the image $x(\{v_1,\dots,v_n\})$, 
where $x_L$ and $x_R$ indicate the left/right extrema. 
We fix $i\ne j\in\{1,\dots,n\}$ such that 
$x(v_i)=x_L$ and $x(v_j)=x_R$ and assign the degree as $|v_i|=|v_j|=0$. 
The degree of the remaining points is set to be one. 
The choice of such $i,j$ is not unique only if 
$v_1=v_n$ and further $x(v_1)=x(v_n)=x_L$ or $x(v_1)=x(v_n)=x_R$. 
Hereafter, by a CC-polygon $\vec{v}$, 
we mean that with a degree attached in the sense above. 
The {\em sign $\sigma(\vec{v})$ of the CC-polygon} $\ti{v}$ 
is then defined by 
\begin{equation*}
 \sigma(\vec{v}):=
 \begin{cases}
   -1 & i < j \\
   +1 & j < i. 
 \end{cases}  
\end{equation*}
(See Figure \ref{fig:polygon-sign}.) 
\begin{figure}[h]
 \begin{minipage}[c]{75mm}{
\begin{center}
 \includegraphics{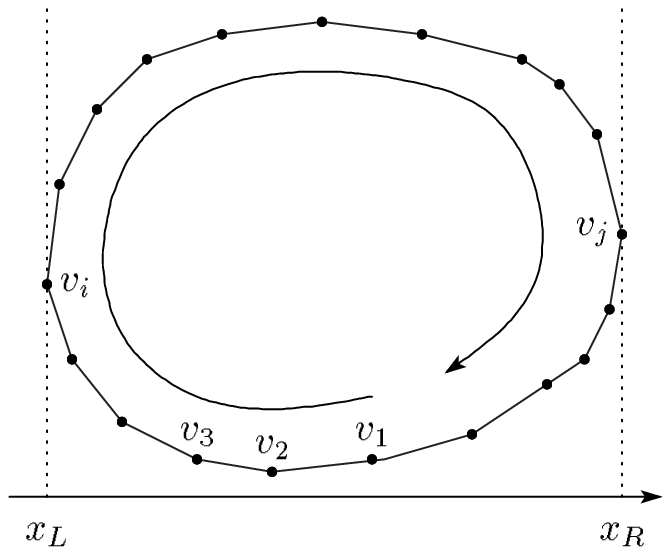}

 (a)
\end{center}}
 \end{minipage}
\quad 
 \begin{minipage}[c]{75mm}{
\begin{center}
 \includegraphics{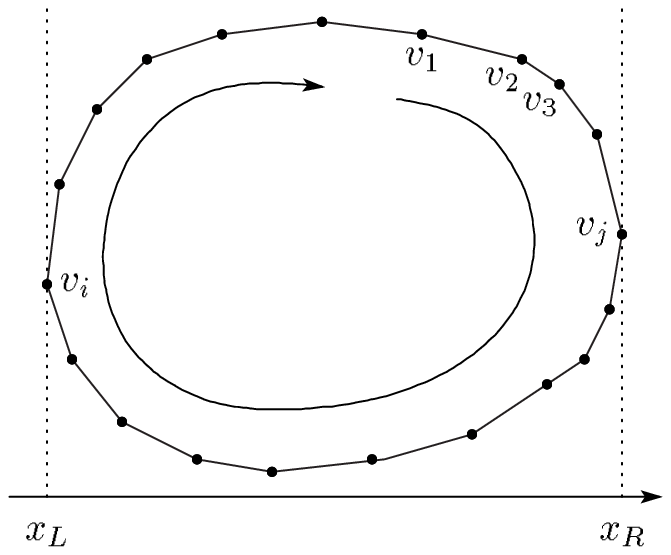}

 (b)
\end{center}}
 \end{minipage}
 \caption{CC-polygons $\vec{v}$ 
with (a): $\sigma(\vec{v})=-1$ and (b): $\sigma(\vec{v})=+1$. }
\label{fig:polygon-sign}
\end{figure}
\end{defn}
For $k\ge 2$, let us define 
degree $(2-k)$ multilinear maps 
$m_k:V_{a_1a_2}\otimes\cdots\otimes V_{a_ka_{k+1}}\to V_{a_1a_{k+1}}$, 
\begin{equation*}
  m_k(w_{a_1a_2},\dots,w_{a_ka_{k+1}})
 \in V_{a_1a_{k+1}}
\end{equation*}
which do not include degree one element in $V_{aa}$ for any $a\in\fF_N$. 
Namely, we consider the case $a_i\ne a_{i+1}$ for any $i=1,\dots,k$, 
which implies that 
$w_{a_ia_{i+1}}$, $i=1,\dots, k$, is spanned by the base $[v_{a_ia_{i+1}}]$. 
Thus, it is enough to determine the multilinear maps 
\begin{equation}\label{determine}
 m_k([v_{a_1a_2}],\dots,[v_{a_ka_{k+1}}]). 
\end{equation}

In the case $a_1\ne a_{k+1}$\ ($k\ge 2$): 
Let us set $\vec{v}:=(v_{a_1a_2},\dots,v_{a_ka_{k+1}},v_{a_{k+1}a_1})$, 
where $v_{a_{i-1}a_{i}}\ne v_{a_{i}a_{i+1}}$ for any $i\in\Z$, $a_i=a_{i+(k+1)}$. 
By degree counting (Lemma \ref{lem:degree}), 
the multilinear maps (\ref{determine}) 
can be nonzero only if 
$\vec{v}$ forms a CC-polygon 
(we shall check this fact in subsection \ref{ssec:polygon-degree}), 
where the degree for the points is attached uniquely as 
$|v_{a_ia_{i+1}}|=|[v_{a_ia_{i+1}}]|$, $i\in\Z$, $a_i=a_{i+(k+1)}$. 
We set the structure constant $c_{a_1\cdots a_{k+1}}\in\R$ of 
\begin{equation*}
 m_k([v_{a_1a_2}],\dots,[v_{a_ka_{k+1}}])
 = c_{a_1\cdots a_{k+1}} [v_{a_1a_{k+1}}]
\end{equation*}
by 
\begin{equation}\label{str-const1}
 c_{a_1\cdots a_{k+1}}:=
 (\sigma(\vec{v}))^k\ e^{-Area(\vec{v})}, 
\end{equation}
where $Area(\vec{v})$ is the area of the CC-polygon $\vec{v}$.

In the case $a_1= a_{k+1}=:a$\ ($k\ge 2$): 
Let us set $\vec{v}:=(v_{a_1a_2},\dots,v_{a_ka_{k+1}})$. 
By degree counting (Lemma \ref{lem:degree}), 
the multilinear maps (\ref{determine}) 
can be nonzero only if 
$\vec{v}$ forms a CC-polygon or a point. 
If $\vec{v}$ forms a CC-polygon, with the degree attached uniquely as 
$|v_{a_ia_{i+1}}|=|[v_{a_ia_{i+1}}]|$, $i=1,\dots,k$, 
we set 
\begin{equation}\label{str-const2}
 \begin{split}
 & m_k([v_{aa_2}],\dots,[v_{a_ka}])
 = c_{a a_2\cdots a_k a}\cdot\ 
 \vartheta_{v_{aa_2}}^{-\sigma(\vec{v})}
 \vartheta_{v_{a_ka}}^{\sigma(\vec{v})}
\in V_{aa}^0\\
 & c_{a_1\cdots a_{k+1}}:=(\sigma(\vec{v}))^k\ e^{-Area(\vec{v})}, 
 \end{split}
\end{equation}
where $\vartheta_{v}^{\pm 1}$ denotes 
\begin{equation*}
 \vartheta_v^{+1}=\vartheta_v,\qquad 
 \vartheta_v^{-1}=1-\vartheta_v. 
\end{equation*} 
If $\vec{v}$ forms a point, the corresponding 
multilinear maps (\ref{determine}) become bilinear one 
$m_2:V_{ab}\otimes V_{ba}\to V_{aa}$ for some $a\ne b\in\fF_N$, 
which we set as 
\begin{equation}\label{str-const3}
 m_2([v_{ab}],[v_{ba}])=\delta_{v_{ab}}\in V_{aa}^1. 
\end{equation}
\begin{thm}
These multilinear maps $m_k$ 
define a unique $A_\infty$-structure in $\cC(\fF_N)$. 
\end{thm}
\begin{pf}
The multilinear maps $m_k$ given above (those which do not 
include elements in $V_{aa}^1$ for some $a\in\fF_N$) is in fact 
compatible with the $A_\infty$-constraint (\ref{Ainfty}); 
this fact can be checked directly. 
On the other hand, 
for any $a\in\fF_N$, elements in $V_{aa}^1=\ti{V}_{aa}^1$ is $m_1$-exact 
in $\ti{V}_{aa}$ and then 
all the $A_\infty$-products including those 
are determined uniquely 
by the $A_\infty$-constraint (\ref{Ainfty}). 
The $A_\infty$-structure on $\ti{V}_{**}$ obtained so 
is in fact closed in $V_{**}$. 

Alternatively, all the compatibility, existence, and the uniqueness of the 
$A_\infty$-structure stated above can also be obtained as a corollary of 
the proof of Proposition (\ref{prop:key}) in subsection \ref{ssec:derive} 
where the $A_\infty$-structure on $\cC(\fF_N)$ 
is derived in the framework of HPT. 
\qed\end{pf}

Note that the $A_\infty$-product (\ref{str-const1}) is just 
the transversal $A_\infty$-product in Condition (ii) 
in Theorem \ref{thm:main}. 
This definition of transversal $A_\infty$-products (\ref{str-const1}) 
{\em agrees in the sign} with that given 
by Polishchuk \cite{Poli:Ainf} 
in the two-tori case. 

As for the $A_\infty$-products of $\sharp(\vartheta)=1$, 
as a consequence, 
the $A_\infty$-products can be nonzero only if 
$k_1=k_2=k_3=0$ for type A and B, 
$k_1=k_3=0$ for type C$_{1}$ with $r=0$ and type C$_{3}$, 
and $k_2=k_4=0$ for type C$_{1}$ with $r=1$ and type C$_{2}$.

 \vspace*{0.3cm}

For applications in the future, 
it should be worth giving the formula for the $A_\infty$-products 
$m_k:V_{a_1a_2}\otimes\cdots\otimes V_{a_ka_{k+1}}\to V_{a_1a_{k+1}}$ 
which do not include elements in $V_{aa}^0$ for any $a$ 
and which may include degree one elements in $V_{aa}^1$ for any $a$ 
of the form $\delta_{v_a}\in V_{aa}^1$ only. 
Namely, for the $A_\infty$-products 
\begin{equation}\label{Ainf-formula}
  m_k(w_{a_1a_2},\dots,w_{a_ka_{k+1}})
 \in V_{a_1a_{k+1}}, \qquad  c_{a_1\cdots a_{k+1}}\in\R, \quad (k\ge 2), 
\end{equation}
we let 
$w_{a_ia_{i+1}}=[v_{a_ia_{i+1}}]\in V_{a_ia_{i+1}}^0$ 
if $t_{a_i}<t_{a_{i+1}}$, 
$w_{a_ia_{i+1}}=[v_{a_ia_{i+1}}]\in V_{a_ia_{i+1}}^1$ 
if $t_{a_i}>t_{a_{i+1}}$ and 
$w_{a_ia_{i+1}}=\delta_{v_{a_i}}\in V_{a_ia_{i+1}}^1$ 
for some $v_{a_i}\in S_{a_i}$ if $a_i=a_{i+1}$, 
where $i=1,\dots,k$. 
Recall that, in this case, any $w_{a_ia_{i+1}}$ is associated with a point 
in $\R^2$. 
If $a_i\ne a_{i+1}$, the associated point is $v_{a_ia_{i+1}}$. 
If $a_i= a_{i+1}$, we denote the associated point again 
by $v_{a_ia_{i+1}}$, where $v_{a_ia_{i+1}}\in S_{a_i}$. 
Then, 
for the case $a_1\ne a_{a_{k+1}}$, we set 
$\vec{v}:=(v_{a_1a_2},\dots,v_{v_{a_ka_{k+1}}},v_{a_{k+1}a_1})$, 
the points associated to elements 
$(w_{a_1a_2},\dots,w_{a_ka_{k+1}},[v_{a_{k+1}a_1}])$. 
For the case $a_1=a_{a_{k+1}}$, we set 
$\vec{v}:=(v_{a_1a_2},\dots,v_{v_{a_ka_{k+1}}})$, 
the points associated to elements $(w_{a_1a_2},\dots,w_{a_ka_{k+1}})$. 
In both cases, the $A_\infty$-product 
$m_k(w_{a_1a_2},\dots,w_{a_ka_{k+1}})$ in eq.(\ref{Ainf-formula}) 
can be zero only if 
$\vec{v}$ forms a CC-polygon or a point. 
Let us describe $\vec{v}$ as the form 
\begin{equation*}
 \vec{v}=(v_1,\dots,v_1,v_2,\dots,v_2,\dots\dots,v_n,\dots,v_n). 
\end{equation*}
Suppose that $\vec{v}$ is a CC-polygon, 
where the degree for each point is given as follows. 
For any $i=1,\dots, n$, if $v_i,\dots,v_i$ includes the point associated to 
a degree zero element, we set $|v_i|=0$ and denote the copy of $v_i$ by 
\begin{equation*}
 (v_i)^{\otimes(d_{i_-},d_{i_+})}:=
 \overbrace{v_i,\dots,v_i}^{d_{i_-}},
 \vv_i,\overbrace{v_i,\dots,v_i}^{d_{i_+}}, 
\end{equation*}
where we attach $\circ$ to the point associated to the degree zero element. 
We put $d_i:=d_{i-}+d_{i_+}$. 
If $v_i,\dots,v_i$ does not include the point associated to 
a degree zero element, we set $|v_i|=1$ 
and denote the $d_i$ copy of $v_i$ by 
\begin{equation*}
 (v_i)^{\otimes d_i}:=\overbrace{v_i,\dots,v_i}^{d_{i}}. 
\end{equation*}
For any point $v_i$ of degree zero or one, 
we call the integer $d_i$ the {\em multiplicity} of $v_i$. 
For any $i=1,\dots, n$, we define 
\begin{equation*}
 D_i= 
 \begin{cases}
  2^{d_i}(d_{i_-})! (d_{i_+})!  & |v_i|=0 \\ 
            (d_i)!  & |v_i|=1. 
 \end{cases} 
\end{equation*}

In the setting above, the $A_\infty$-products (\ref{Ainf-formula}) 
is determined as follows. 

In the case $a_1\ne a_{k+1}$\ ($k\ge 2$): 
The $A_\infty$-product (\ref{Ainf-formula}) is nonzero 
only if $\vec{v}$ forms a CC-polygon (\ie, is zero if $\vec{v}$ is a point). 
Then, the structure constant $c_{a_1\cdots a_{k+1}}\in\R$ of 
\begin{equation*}
 m_k(w_{a_1a_2},\dots,w_{a_ka_{k+1}})
 = c_{a_1\cdots a_{k+1}} [v_{a_1a_{k+1}}]
\end{equation*}
is given by 
\begin{equation}\label{str-const1'}
 c_{a_1\cdots a_{k+1}}:=
 \frac{(\sigma(\vec{v}))^k}{D_1\cdots D_n}\ e^{-Area(\vec{v})}. 
\end{equation}

In the case $a_1= a_{k+1}=:a$\ ($k\ge 3$): 
The $A_\infty$-product (\ref{Ainf-formula}) can be nonzero 
only if $\vec{v}$ forms a CC-polygon or a point. 
If $\vec{v}$ forms a CC-polygon, it is given by 
\begin{equation}\label{str-const2'}
 \begin{split}
 & m_k(w_{aa_2},\dots,w_{a_ka})
 = c_{a a_2\cdots a_k a}\cdot\alpha_1(\vec{v})\alpha_n(\vec{v}) \in V_{aa}^0\\
 & c_{a_1\cdots a_{k+1}}:=
 \frac{(\sigma(\vec{v}))^k}{D_1\cdots D_n}
\ e^{-Area(\vec{v})}, 
 \end{split}
\end{equation}
where $\alpha_1(\vec{v}), \alpha_n(\vec{v})\in V_{aa}^0$ are defined by 
\begin{equation*}
 \alpha_1(\vec{v})=
 \begin{cases}
  \vartheta_{v_1}^{-\sigma(\vec{v})}
  (-\sigma(\vec{v})2(\vartheta_{v_1}-\ov{2}))^{d_1} & |v_1|=0, \\
  (\vartheta_{v_1}^{-\sigma(\vec{v})} )^{d_1} & |v_1|=1, 
 \end{cases}
\end{equation*}
\begin{equation*}
 \alpha_n(\vec{v})=
 \begin{cases}
  \vartheta_{v_n}^{\sigma(\vec{v})}
  (\sigma(\vec{v})2(\vartheta_{v_n}-\ov{2}))^{d_n} & |v_n|=0, \\
  (\vartheta_{v_n}^{\sigma(\vec{v})} )^{d_n} & |v_n|=1. 
 \end{cases}
\end{equation*}
If $\vec{v}$ forms a point, 
nonzero $A_\infty$-products are only the followings: 
\begin{equation*}
 \begin{split}
 & m_{2+d_-+d_+}(\overbrace{\delta_{v_{ab}},\dots,\delta_{v_{ab}}}^{d_-},[v_{ab}],
 \overbrace{\delta_{v_{ab}},\dots,\delta_{v_{ab}}}^{d_+},[v_{ba}])
 =\frac{\left(\ov{2}-\vartheta_{v_{ab}}\right)^{d_-+d_+}}{(d_-)!(d_+)!}
 \cdot\delta_{v_{ab}}\in V_{aa}^1, \\
 & m_{2+d_-+d_+}([v_{ba}],
\overbrace{\delta_{v_{ab}},\dots,\delta_{v_{ab}}}^{d_-},[v_{ab}],
 \overbrace{\delta_{v_{ab}},\dots,\delta_{v_{ab}}}^{d_+})
 =\frac{\left(\ov{2}-\vartheta_{v_{ab}}\right)^{d_-+d_+}}{(d_-)!(d_+)!}
 \cdot\delta_{v_{ab}}\in V_{bb}^1 
 \end{split}
\end{equation*}
for any $a, b\in\fF_N$ such that $t_a<t_b$ 
and $d_-,d_+\in\Z_{\ge 0}$. 
\begin{thm}\label{thm:equiv}
For any given $\fF_N$ and $\fF_N'$, 
the two $A_\infty$-categories $\cC(\fF_N)$ and $\cC'(\fF'_N)$ 
are homotopy equivalent. 
\end{thm}
We shall prove this after the proof of the main theorem (Theorem \ref{thm:main}) 
in subsection \ref{ssec:proof}. 
This result seems reasonable from the viewpoint of symplectic geometry, 
since the second cohomology of $\R^2$ is trivial, that is, 
any symplectic form on $\R^2$ is exact. 
Extending this construction of an $A_\infty$-structure to tori, 
the $A_\infty$-categories with different configurations of lines 
are not homotopy equivalent in general 
even if the number of the lines (objects) is the same. 
For instance, for a two-torus case, 
a geodesic cycle in a torus $T^2$ is described by a $\Z$-copy of 
lines in the covering space $\R^2$ and 
the dimension of the space $\Hom_{Fuk(T^2)}(a,b)$ of morphisms 
for two transversal geodesic cycles $a$ and $b$ in $T^2$ 
is the number of the intersection points of $a$ and $b$ in $T^2$, 
which will change if we change the slops of $a$ and $b$ even if we keep 
the ordering (cf. \cite{PZ,hk:foliation, PoSc,hk:nchms}).

 \section{Interpretations and Examples}
\label{sec:interpret}

In this section, we explain more on the relation of polygons and 
the $A_\infty$-products in the Fukaya $A_\infty$-category $\cC(\fF_N)$ 
mainly for the transversal part. 
The relation between polygons and the degrees 
of the the intersection points 
is explained in subsection \ref{ssec:polygon-degree}. 
The realization of the $A_\infty$-constraints in terms of polygons 
is given in subsection \ref{ssec:polygon-relation}. 
Also, the necessity of nontrivial non-transversal products is observed 
in an example in subsection \ref{ssec:why}.

 \subsection{CC-polygon and the degree}
\label{ssec:polygon-degree}

In the previous subsection we stated that 
by degree counting (Lemma \ref{lem:degree}) 
the transversal $A_\infty$-products 
$m_n([v_{a_1a_2}],\dots,[v_{a_na_{n+1}}])$ 
in eq. (\ref{str-const1}) can be nonzero 
only if the corresponding sequence 
$\vec{v}=(v_{a_1a_2},\dots,v_{a_na_{n+1}})$ of points 
forms a CC-polygon. Let us check this fact. 
\begin{figure}[h]
 \includegraphics{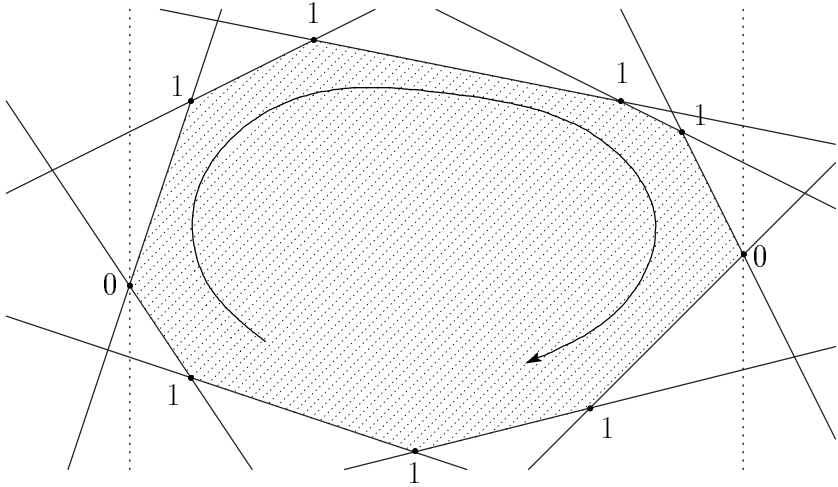}
 \caption{A CC-polygon $\vec{v}$ 
with degree ($0$ or $1$) of the intersection points. }\label{fig:polygon}
\end{figure}
If we go around the CC $(n+1)$-gon $\vec{v}$ in the clockwise direction 
and count the degree $r$ (zero or one) 
of each point $v_{a_ia_{i+1}}$, 
we always have two points of degree zero and $(n+1)-2$ points 
of degree one as in Figure \ref{fig:polygon}. 
Thus, we have 
\begin{equation}\label{degree-transversal}
 \sum_{i=1}^{n+1} |v_{a_ia_{i+1}}|=
 \sum_{i=1}^{n+1} |[v_{a_ia_{i+1}}]|=(n+1)-2 
\end{equation}
for the CC $(n+1)$-gon $\vec{v}$, where $v_{a_{n+1}a_{n+2}}:=v_{a_{n+1}a_1}$. 
One can also see that the equation above holds true only if $\vec{v}$ 
forms a CC-polygon. 

On the other hand, for the transversal $A_\infty$-product 
$m_n([v_{a_1a_2}],\dots,[v_{a_na_{n+1}}])$, 
Lemma \ref{lem:degree} implies that one has only two elements of degree zero 
in $\{[v_{a_1a_2}],\dots,[v_{a_na_{n+1}}],[v_{a_{n+1}a_1}]\}$. 
This exactly implies the identity (\ref{degree-transversal}) 
since the degree of the remaining elements is one. 

To make sure, let us check Lemma \ref{lem:degree} 
in this transversal situation. 
Let us assume that $|[v_{a_{n+1}a_1}]|=r$, where $r$ is equal to zero or one. 
Then, one has 
\begin{equation*}
 |m_n([v_{a_1a_2}],\dots,[v_{a_na_{n+1}}])|
 =\sum_{i=1}^{n} |[v_{a_ia_{i+1}}]| -r + (2-n) 
\end{equation*}
since the degree of $m_n$ is $(2-n)$. 
Here, $\vec{v}$ is a CC-polygon if and only if 
the identity (\ref{degree-transversal}) holds, 
where the right hand side of the equation above 
turns out to be $(n+1)-2+(2-n)-r=1-r= |[v_{a_1a_n}]|$. 

Thus, one can indeed define 
nonzero transversal $A_\infty$-product $m_n$ only when 
the corresponding $\vec{v}$ forms a CC $(n+1)$-gon. 
This fact was obtained by counting the degrees of the corresponding points, 
which are related to their Maslov indices (see \cite{fukayaAinfty}. )

 \subsection{$A_\infty$-constraint and polygons}
\label{ssec:polygon-relation}

In the rest of this section, we denote $v_{ab}:=[v_{ab}]$ since 
it does not cause any confusion. 

The $A_\infty$-constraints for transversal $A_\infty$-products 
has a geometric interpretation in terms of a clockwise polygon which 
has one nonconvex point (=vertex of the polygon). 
There exist two ways to divide the polygon into two convex polygons. 
The corresponding terms then appear with opposite signs and 
cancel each other in the $A_\infty$-constraint. 
\begin{figure}[h]
 \includegraphics{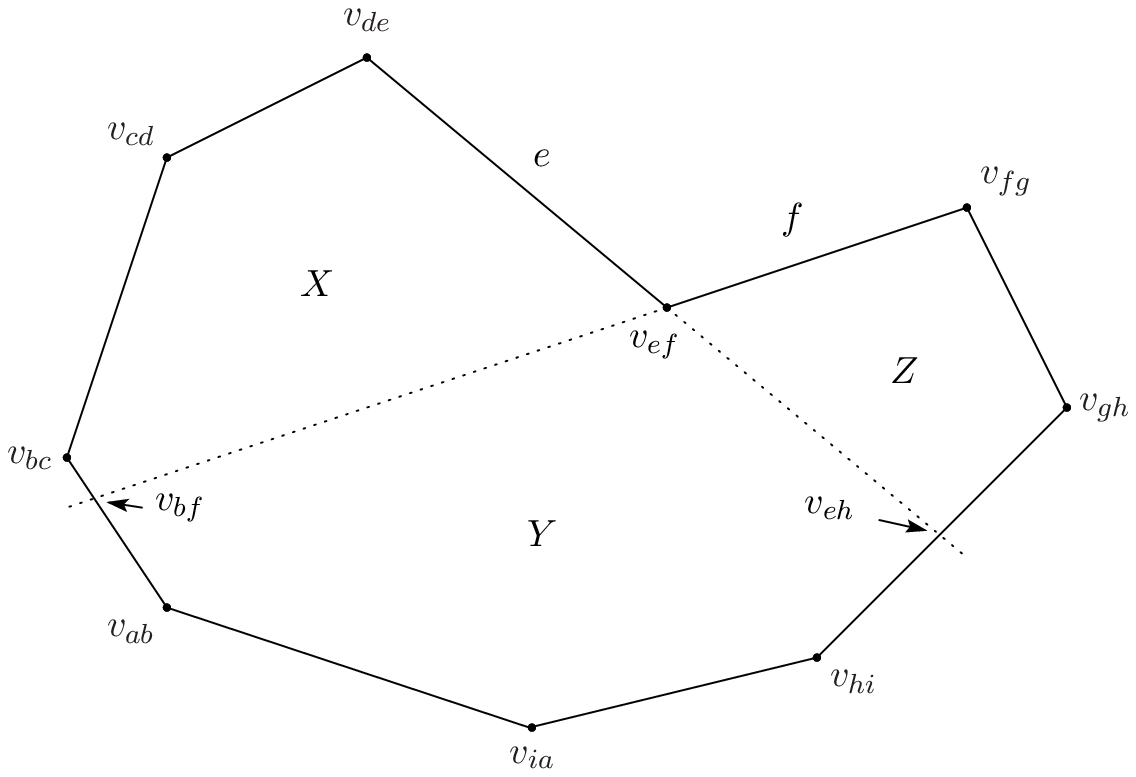}
 \caption{A clockwise polygon 
which has one nonconvex point. }\label{fig:polygon-rel}
\end{figure}
For example, in Figure \ref{fig:polygon-rel}, we have the following 
(intersection) points with their degrees assigned: 
\begin{equation*}
 \begin{array}{cccccccc}
  v_{ab} &v_{bc}&v_{cd} &v_{de} &v_{ef} &v_{fg} &v_{gh} &v_{hi} \\
    1    & 0    & 1     & 1     & 0     & 1     & 0     & 1
 \end{array}\ .
\end{equation*}
Corresponding to the way of dividing the area $X+Y+Z$ into 
(i) $X+(Y+Z)$ or (ii) $(X+Y)+Z$, we have the following composition 
of transversal $A_\infty$-products: 
\begin{equation*}
 \begin{split}
 \mbox{(i)} &\pm  (v_{ab}(v_{bc}v_{cd}v_{de}v_{ef})v_{fg}v_{gh}v_{hi}), \\
 \mbox{(ii)}&\pm  (v_{ab}v_{bc}v_{cd}v_{de}(v_{ef}v_{fg}v_{gh})v_{hi}),
 \end{split}
\end{equation*}
where $(v_{bc}v_{cd}v_{de}v_{ef})$ indicates 
$m_4(v_{bc},v_{cd},v_{de},v_{ef})$ and so on. 
There does not exist any other composition of 
$A_\infty$-products since a transversal $A_\infty$-product can be nonzero 
only if the corresponding polygon forms a CC-polygon. 
According to the definition, one obtains 
\begin{equation*}
 \begin{split}
 &m_4(v_{bc},v_{cd},v_{de},v_{ef})=e^{-X}v_{bf}\ ,\\
 &m_5(v_{ab},v_{bf},v_{fg},v_{gh},v_{hi})=-e^{-(Y+Z)}v_{ai}\ ,\\
 &m_3(v_{ef},v_{fg},v_{gh})=-e^{-Z}v_{eh}\ ,\\
 &m_6(v_{ab},v_{bc},v_{cd},v_{de},v_{eh},v_{hi})=e^{-(X+Y)}v_{ai}\ .
 \end{split}
\end{equation*}
Combining the first two equations leads to 
\begin{equation*}
 m_5(v_{ab},m_4(v_{bc},v_{cd},v_{de},v_{ef}),v_{fg},v_{gh},v_{hi})
 =-e^{-X-(Y+Z)}v_{ai}\ ,
\end{equation*}
and combining the last two gives 
\begin{equation*}
 m_6(v_{ab},v_{bc},v_{cd},v_{de},m_3(v_{ef},v_{fg},v_{gh}),v_{hi})
 =-e^{-(X+Y)-Z}v_{ai}\ .
\end{equation*}
Thus, we obtain 
\begin{equation*}
0=
m_5(v_{ab},m_4(v_{bc},v_{cd},v_{de},v_{ef}),v_{fg},v_{gh},v_{hi})
-m_6(v_{ab},v_{bc},v_{cd},v_{de},m_3(v_{ef},v_{fg},v_{gh}),v_{hi})\ ,
\end{equation*}
which is just the $A_\infty$-constraints (\ref{Ainfty}) 
on $(v_{ab},v_{bc},v_{cd},v_{de},v_{ef},v_{fg},v_{gh},v_{hi})$. 

 \subsection{Why can we not avoid non-transversal products ? }
\label{ssec:why}

Using Figure \ref{fig:polygon-rel}, we show that we can not 
avoid non-transversal $A_\infty$-products, \ie, 
we can not define an $A_\infty$-structure for the Fukaya category 
such that all non-transversal $A_\infty$-products are zero. 

Consider the sequence 
$(v_{ab},v_{bf},v_{fe},v_{ef},v_{fg},v_{gh},v_{hi})$ of 
elements and the corresponding $A_\infty$-products. 
There exists a composition 
$m_5(v_{ab},v_{bf},v_{fe},m_3(v_{ef},v_{fg},v_{gh}),v_{hi})
=e^{-(Y+Z)}v_{ai}$ 
of two transversal $A_\infty$-products. The $A_\infty$-constraint (\ref{Ainfty}) 
then implies that this composition cancels with other terms. 
However, there does not exist any more composition 
of two non-zero transversal $A_\infty$-products 
on the sequence $(v_{ab},v_{bf},v_{fe},v_{ef},v_{fg},v_{gh},v_{hi})$. 
This shows the necessity of nonzero non-transversal $A_\infty$-products. 
For the case of the $A_\infty$-category $\cC(\fF_N)$, 
one has $m_2(v_{fe},v_{ef})=\delta_{v_{fe}}\in V_{ff}^1$ and 
the $A_\infty$-constraint on the sequence 
$(v_{ab},v_{bf},v_{fe},v_{ef},v_{fg},v_{gh},v_{hi})$ is 
\begin{equation*}
 0= m_5(v_{ab},v_{bf},v_{fe},m_3(v_{ef},v_{fg},v_{gh}),v_{hi})
 - m_6(v_{ab},v_{bf},m_2(v_{fe},v_{ef}),v_{fg},v_{gh},v_{hi}), 
\end{equation*}
where 
$m_6(v_{ab},v_{bf},m_2(v_{fe},v_{ef}),v_{fg},v_{gh},v_{hi})
=m_6(v_{ab},v_{bf},\delta_{v_{fe}},v_{fg},v_{gh},v_{hi})
=e^{-(Y+Z)}v_{ai}$.

 \section{Proof of Theorem \ref{thm:main}}
\label{sec:proof}

 \subsection{The outline of the proof}
\label{ssec:proof}

We define two DG categories $\cC_{DR}'(\fF_N)$ and $\ti{\cC}_{DR}(\fF_N)$, 
which are DG-categorical extensions of 
DG-algebras $A_S(\R)$ and $A_S(\Omega(\R))$, respectively. 

Let $S_{all}$ be the set of all intersection points $v_{ab}$ 
for all $a\ne b\in\fF_N$ equipped with a map 
$x: S_{all}\to\R$. Here, note that $v_{ab}=v_{ba}\in S_{all}$, 
$v_{ab}=v_{ac}$ if and only if $b=c\in\fF_N$, and 
in general $x(v)=x(v')$ possibly holds for $v\ne v'\in S_{all}$. 
\begin{defn}[$\ti{\cC}_{DR}(\fF_N)$]\label{defn:tiDG}
The set of objects is taken to be the same as that of $\cC_{DR}(\fF_N)$: 
\begin{equation*}
 \Ob(\ti{\cC}_{DR}(\fF_N)):=\Ob(\cC_{DR}(\fF_N))=\fF_N .  
\end{equation*}
For any $a,b\in\fF_N$, we set the space of morphisms by 
\begin{equation*}
 \Hom_{\ti{\cC}_{DR}(\fF_N)}(a,b)=\ti{\Omega}_{ab}:=A_{ S_{all}}(\Omega(\R))  
\end{equation*}
as a graded vector space of degree zero and one. 
For $a,b,c\in\fF$, the composition 
$m:\ti{\Omega}_{ab}\otimes\ti{\Omega}_{bc}\to\ti{\Omega}_{ac}$ 
is defined as the product in $A_{ S_{all}}(\Omega(\R))$. 
For $a,b\in\fF_N$, 
the differential $d_{ab}:\ti{\Omega}_{ab}\to\ti{\Omega}_{ab}$ 
is given by 
\begin{equation*}
 d_{ab}=d-df_{ab}\wedge, 
\end{equation*}
where $d:A_{ S_{all}}(\Omega(\R))\to A_{ S_{all}}(\Omega(\R))$ 
is the differential of $A_{ S_{all}}(\Omega(\R))$, 
and $\wedge$ is the graded commutative product 
in $A_{ S_{all}}(\Omega(\R))$. 
\end{defn}
\begin{defn}[$\cC_{DR}'(\fF_N)$]\label{defn:DG'}
The set of objects is the same as that of $\cC_{DR}(\fF_N)$: 
\begin{equation*}
 \Ob(\cC'_{DR}(\fF_N)):=\fF_N .  
\end{equation*}
For any $a,b\in\fF_N$, we set the space of morphisms by 
\begin{equation*}
 \Hom_{\cC'_{DR}(\fF_N)}(a,b)
 =\Omega'_{ab}:=\{e^{f_{ab}}\cdot\alpha\in A_{S_{all}}(\Omega(\R))\ |
   \ \alpha\in \ti{A}_{S_{all}}(\R)\} 
\end{equation*}
as a graded vector space of degree zero and one. 
For $a,b,c\in\fF$, the composition 
$m:\Omega'_{ab}\otimes\Omega'_{bc}\to\Omega'_{ac}$ 
is defined as the product in $A_{ S_{all}}(\Omega(\R))$. 
For $a,b\in\fF_N$, 
the differential $d_{ab}:\Omega'_{ab}\to\Omega'_{ab}$ 
is the same as that in $\cC_{DR}(\fF_N)$ or $\ti{\cC}_{DR}(\fF_N)$: 
\begin{equation*}
 d_{ab}=d-df_{ab}\wedge. 
\end{equation*} 
\end{defn}
Clearly, $\ti{\cC}_{DR}(\fF_N)$ forms a DG category, 
and the DG subcategory $\cC_{DR}'(\fF_N)\subset\ti{\cC}_{DR}(\fF_N)$ 
is well-defined. 

The following can be showed just in the same way as Lemma \ref{lem:CDG}: 
\begin{lem}\label{lem:DG-cat}
The inclusions 
\begin{equation*}
 \iota: \cC_{DR}'(\fF_N)\to\ti{\cC}_{DR}(\fF_N),\qquad 
 \iota: \cC_{DR}(\fF_N)\to\ti{\cC}_{DR}(\fF_N),
\end{equation*}
induce homotopy equivalences 
$\cC_{DR}(\fF_N)\simeq\ti{\cC}_{DR}(\fF_N)$ and 
$\cC'_{DR}(\fF_N)\simeq\ti{\cC}_{DR}(\fF_N)$ 
as $A_\infty$-categories. 
\qed\end{lem}
On the other hand, one has the following: 
\begin{prop}\label{prop:key}
There exists an $A_\infty$-functor $\cG:\cC(\fF_N)\to\cC'_{DR}(\fF_N)$ 
which induces a homotopy equivalence 
\begin{equation*}
 \cC(\fF_N)\simeq\cC'_{DR}(\fF_N). 
\end{equation*}
\end{prop}
We shall show this in the next subsection, where HPT is applied 
to derive the $A_\infty$-structure of $\cC(\fF_N)$. 

Then, we obtain the following homotopy equivalences 
\begin{equation*}
 \cC(\fF_N)\overset{\cG}{\to} \cC_{DR}'(\fF_N)
 \overset{\iota}{\to}\ti{\cC}_{DR}(\fF_N) 
 \overset{\iota}{\law}\cC_{DR}(\fF_N), 
\end{equation*}
which give a proof of Theorem \ref{thm:main}. \qed

\begin{pf} {\bf of Theorem \ref{thm:equiv}.}\ \ 
By Theorem \ref{thm:main} $\cC(\fF_N)\simeq\cC_{DR}(\fF_N)$ 
and $\cC(\fF_N')\simeq \cC_{DR}(\fF_N')$. 
Also, by Lemma \ref{lem:DG-cat} we obtained the equivalence 
$\cC_{DR}(\fF_N)\simeq\cC'_{DR}(\fF_N)$ 
and $\cC_{DR}(\fF_N')\simeq \cC'_{DR}(\fF_N')$. 

Thus, one may show the homotopy equivalence 
$\cC'_{DR}(\fF_N)\simeq \cC'_{DR}(\fF_N')$. 
These two categories are in fact isomorphic to each other. 
Let us denote the objects by $\fF_N=\{a,b,\dots\}$ and 
$\fF_N'=\{a',b',\dots\}$, where we can assume $t_a<t_b\cdots$ 
and $t_{a'}<t_{b'}\cdots$ without lose of generality. 
The functor between the objects is given by 
$a\mapsto a'$ for any $a\in\fF_N$, 
and the functor between the space of morphisms is given by 
\begin{equation*}
\renewcommand{\arraystretch}{1.5}
 \begin{array}{ccc}
  \Hom_{\cC'_{DR}(\fF_N)}(a,b) & \to &\Hom_{\cC'_{DR}(\fF_N')}(a',b') \\
     \omega & \mapsto & e^{f_{a'b'}-f_{ab}}\omega 
 \end{array}
\end{equation*}
for any $a,b\in\fF_N$. 
\qed\end{pf}

 \subsection{Deriving the $A_\infty$-category $\cC(\fF_N)$}
\label{ssec:derive}

Now, we shall show Proposition \ref{prop:key} stating 
the homotopy equivalence $\cC(\fF_N)\simeq \cC_{DR}'(\fF_N)$. 
We apply HPT (Theorem \ref{thm:HPT-cat}) to $\cC_{DR}'(\fF_N)$. 
In order to do so, for any $a,b\in\fF_N$, we first define 
homotopy $h_{ab}:{\Omega'}_{ab}^1\to {\Omega'}_{ab}^0$, 
${\Omega'}_{ab}:=\Hom_{\cC_{DR}'(\fF_N)}(a,b)$, 
so that $P_{ab}:{\Omega'}_{ab}^r\to{\Omega'}_{ab}^r$, $r=0,1$, defined by 
the $d_{ab}h_{ab}+h_{ab}d_{ab}=\Id_{\Omega'_{ab}}-P_{ab}$ 
gives a projection on ${\Omega'}_{ab}$. 

For any $a\in\fF_N$, we set $h_{aa}=0$ and then $P_{aa}=\Id$. 

For all $a\ne b\in\fF_N$, define 
the homotopy $h_{ab}: {\Omega'}_{ab}^1\to {\Omega'}_{ab}^0$ 
and the projection $P_{ab}:{\Omega'}_{ab}^r\to {\Omega'}_{ab}^r$ 
as follows. 
For the base $d(\vartheta_v)^n\in {\Omega'}_{ab}^1$, 
\begin{equation}\label{h-ab}
 h_{ab}(d(\vartheta_v)^n):=e^{f_{ab}-f_{ab}(x(v))}((\vartheta_v)^n-c) 
\end{equation}
for the case $t_a<t_b$, where 
$c=0$ if $x(v_{ab})<x(v)$, $c=1/2^n$ if $x(v_{ab})=x(v)$, 
and $c=1$ if $x(v)<x(v_{ab})$. 
On the other hand, if $t_a>t_b$, 
\begin{equation}\label{h-ba}
 h_{ab}(d(\vartheta_v)^n)
=e^{f_{ab}-f_{ab}(x(v))}((\vartheta_v)^n-\vartheta_{v_{ab}}) 
\end{equation}
for any $n\ge 1$ and $v\in S_{all}$. 

For $t_a<t_b$, 
the projection $P_{ab}: {\Omega'}_{ab}^0\to {\Omega'}_{ab}^0$ is 
\begin{equation*}
 P_{ab}(e^{f_{ab}}\cdot (\vartheta_v)^n)
 =\begin{cases}
   e^{f_{ab}} & x(v)< x(v_{ab}) \\
   \ov{2^n} e^{f_{ab}} & x(v)= x(v_{ab}) \\
    0 & x(v_{ab}) < x(v) 
 \end{cases}
\end{equation*}
for any $n\ge 1$ and $v\in S_{all}$, $P_{ab}(e^{f_{ab}})=e^{f_{ab}}$, 
and $P_{ab}=0$ for $P_{ab}: {\Omega'}_{ab}^1\to {\Omega'}_{ab}^1$.  

For $t_a>t_b$, 
the projection $P_{ab}: {\Omega'}_{ab}^1\to {\Omega'}_{ab}^1$ is 
\begin{equation*}
 P_{ab}(d(\vartheta_v)^n)
 = \delta_{v_{ab}}
\end{equation*}
for any $n\ge 1$ and $v\in S_{all}$ 
and $P_{ab}=0$ for $P_{ab}: {\Omega'}_{ab}^0\to {\Omega'}_{ab}^0$.  

Then, for any $a,b\in\fF_N$, one has the identity 
\begin{equation}\label{HKdecomp-apply}
 d_{ab}h_{ab}+h_{ab}d_{ab}=\Id-P_{ab}
\end{equation}
on $\Omega'_{ab}$. 
We denote the base of $P_{ab}{\Omega'}_{ab}^r$ by 
$\be_{ab}:=e^{f_{ab}-f_{ab}(x(v_{ab}))}$ for $r=0$ 
and $\be_{ab}:=\delta_{v_{ab}}$ for $r=1$. 

Now, applying HPT (Theorem \ref{thm:HPT-cat}) 
to $\cC'_{DR}(\fF_N)$ with the identity (\ref{HKdecomp-apply}) 
leads to an $A_\infty$-category $\cC'(\fF_N)$ 
with homotopy equivalence $\cC'(\fF_N)\overset{\sim}{\to}\cC'_{DR}(\fF_N)$. 
Here, the set of the objects is $\Ob(\cC'(\fF_N))=\fF_N$. 
The space of morphisms is defined so that 
\begin{equation}\label{iota_ab}
 \iota_{ab}: \Hom_{\cC'(\fF_N)}(a,b)\to\Hom_{\cC'_{DR}(\fF_N)}(a,b)
\end{equation}
gives the embedding to 
$P_{ab}\Hom_{\cC'_{DR}(\fF_N)}(a,b)\subseteq\Hom_{\cC'_{DR}(\fF_N)}(a,b)$ 
for any $a,b\in\fF_N$, which turns out to be 
$\Hom_{\cC'(\fF_N)}(a,b)=\Hom_{\cC(\fF_N)}(a,b)=V_{ab}$ 
for any $a\ne b\in\fF_N$ but 
$\Hom_{\cC'(\fF_N)}(a,a)=\Hom_{\cC'_{DR}(\fF_N)}(a,a)
\supset\Hom_{\cC(\fF_N)}(a,a)=V_{aa}$ 
for any $a\in\fF_N$. 
For $a\ne b\in\fF_N$, 
we identify the base $\be_{ab}\in P_{ab}\Hom_{\cC'_{DR}(\fF_N)}(a,b)$ 
with the base $[v_{ab}]\in\Hom_{\cC'(\fF_N)}(a,b)=V_{ab}$ by 
the embedding $\iota_{ab}$ in eq.(\ref{iota_ab}) as 
$\be_{ab}=\iota_{ab}[v_{ab}]$. 
Then, the $A_\infty$-structure on $\cC'(\fF_N)$ is closed 
in the subspace $\Hom_{\cC(\fF_N)}(a,b)\subseteq\Hom_{\cC'(\fF_N)}(a,b)$, 
which gives the $A_\infty$-structure on $\cC(\fF_N)$. 
Clearly, the inclusion $\cC(\fF_N)\to\cC'(\fF_N)$ gives an $A_\infty$-functor, 
and in fact gives homotopy equivalence 
since the inclusion is a quasi-isomorphism of chain complexes 
$\Hom_{\cC(\fF_N)}(a,b)\to\Hom_{\cC'(\fF_N)}(a,b)$. \qed

Lastly, we end with deriving some examples of the 
$A_\infty$-products (\ref{str-const1'}) (\ref{str-const2'}) 
of $\cC(\fF_N)$ associated with CC-polygons. 

Consider an $A_\infty$-product 
\begin{equation*}
 m_k(w_{a_1a_2},\dots,w_{a_ka_{k+1}}) 
\end{equation*}
which is associated with a CC-polygon $\vec{v}$ 
in the sense in eq.(\ref{Ainf-formula}) and the descriptions below. 
We describe the CC-polygon as 
\begin{equation*}
\vec{v}=((v_1)^{\otimes d_1},(v_2)^{\otimes d_2},\dots, 
(\vv_i)^{\otimes (d_{i_-},d_{i_+})},\dots, 
(\vv_j)^{\otimes (d_{j_-},d_{j_+})},\dots,(v_n)^{\otimes d_n}) 
\end{equation*}
with the map $x:\{v_1,\dots,v_n\}\to\R$ together. 
As above, we attach $\circ$ on degree zero points $v_i$ and $v_j$ 
to distinguish them from other points. 
If $d=1$ for $(v)^{\otimes d}$, we simply denote it as 
$(v)^{\otimes 1}=v$. 
Similarly, if $(d_-,d_+)=(0,0)$ for $(\vv)^{\otimes (d_-,d_+)}$, 
we denote it as $(\vv)^{\otimes (0,0)}=\vv$. 
We shall derive these $A_\infty$-products 
by applying HPT (Theorem \ref{thm:HPT-cat}) 
to the DG category $\cC'_{DR}(\fF_N)$. 
We denote the product in $\cC'_{DR}(\fF_N)$ by $m$. 
To simplify the formula, in these examples, we identify 
$[v_{ab}]\in V_{ab}$ with $\be_{ab}\in P_{ab}\Omega'_{ab}$ for 
$a\ne b\in\fF_N$, and then the surjection 
$\pi_{ab}:\Omega'_{ab}\to V_{ab}$ with $P_{ab}$.

Let us start from deriving a transversal $A_\infty$-product with an example. 
One can see how the area of the corresponding CC-polygon appears, 
where the correspondence of the CC-polygon and 
a planar tree (a Feynman graph) is a key point. 
The way of determining the sign shall be explained 
in the end of this subsection. 
Therefore, in the examples below, we do not care about the sign 
and denote it simply by $\pm$. 
\begin{exmp}\label{exmp:1}
$\vec{v}
=(\vv_{ab},v_{bc},\vv_{cd},v_{da})$, 
$x(v_{ab})<x(v_{bc})<x(v_{da})<x(v_{cd})$. 
The HPT implies 
\begin{equation}\label{tree-transversal}
m_3(\be_{ab},\be_{bc},\be_{cd})\ \ = \ 
\begin{minipage}[c]{36mm}{\includegraphics{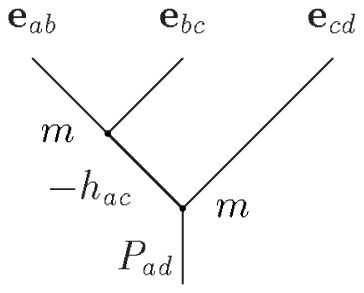}}
\end{minipage}
\ +\ 
\begin{minipage}[c]{36mm}{\includegraphics{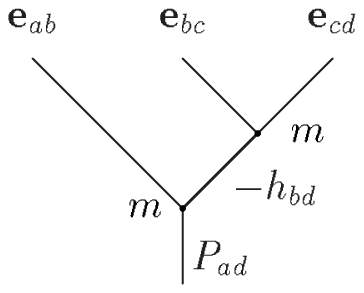}}
\end{minipage}. 
\end{equation}
On the other hand, one can associate a planar tree graph 
to the CC-polygon $\vec{v}$ as follows. 
First, connect two points $v_{ab}, v_{cd}$ of $\vec{v}$ 
of degree zero with an interval. 
For each point of $\vec{v}$ of degree one, 
draw an interval (external edge = leaf), perpendicular to the $x$-axis, 
starting from the point and ending on the interval $(v_{ab}v_{cd})$. 
Choosing the interval starting from the point $v_{da}$ as 
the root edge, one obtains a planar rooted tree 
as in Figure \ref{fig:m-transversal}. 
\begin{figure}[h]
 \includegraphics{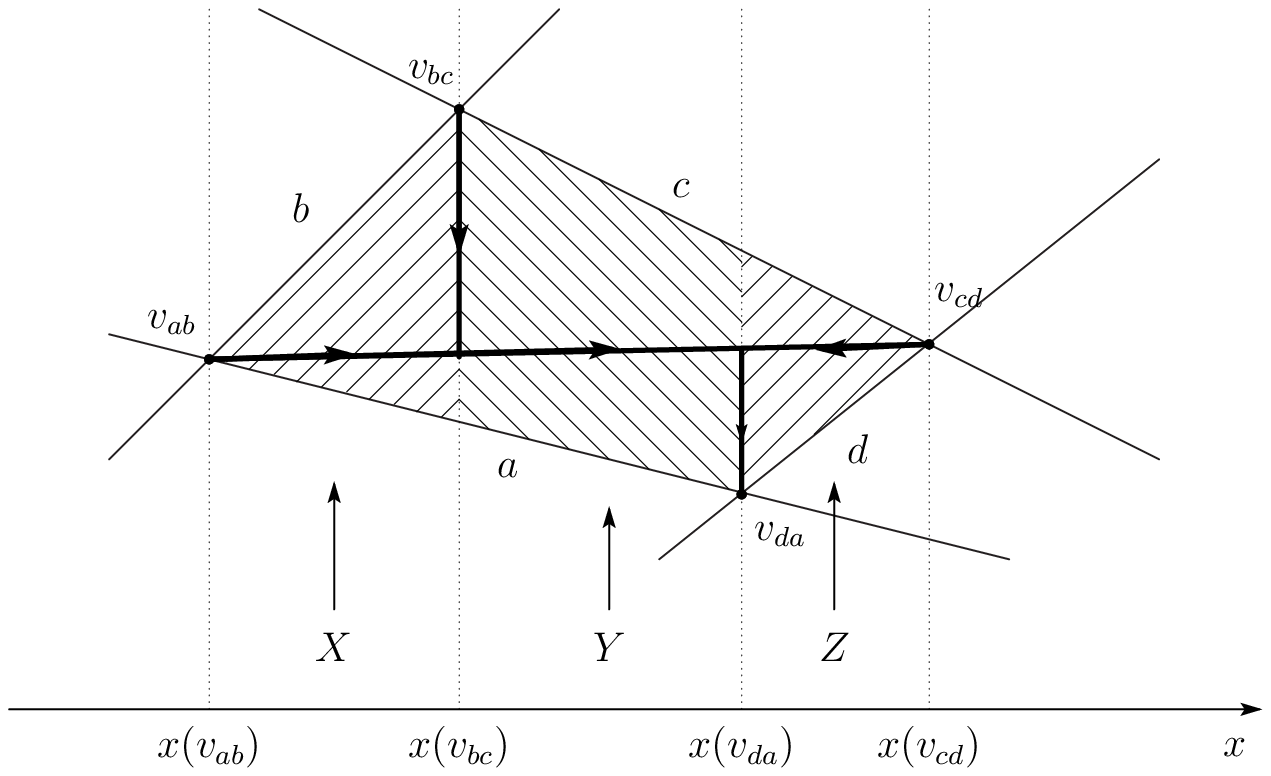}
 \caption{CC-polygon $\vec{v}
=(\vv_{ab},v_{bc},\vv_{cd},v_{da})$. }\label{fig:m-transversal}
\end{figure}
One can see that the resulting planar rooted tree corresponds to 
the one in the first term of the right hand side 
of eq.(\ref{tree-transversal}). 
We shall show that the second term of the right hand side 
of eq.(\ref{tree-transversal}) in fact vanishes and the first term 
derives the area $Area(\vec{v})$. 
Let us calculate the first term. 
As in Figure \ref{fig:m-transversal}, 
we divide the CC-polygon $\vec{v}$ into three by the lines through 
$v_{bc}$ and $v_{da}$ both of which are perpendicular to the $x$-axis. 
The areas between $x(v_{ab})$ and $x(v_{bc})$, $x(v_{bc})$ and $x(v_{da})$, 
$x(v_{da})$ and $x(v_{cd})$ are denoted $X$, $Y$, $Z$, respectively. 
First, one gets $m(\be_{ab},\be_{bc})=\pm e^{-X}\delta_{v_{bc}}$. 
We know $h_{ac}\delta_{v_{bc}}=\pm e^{f_{ac}-f_{ac}(x(v_{bc}))}\cdot\vartheta_{v_{bc}}$. 
Then, $P_{ad} m(-h_{ac}\delta_{v_{bc}},\be_{cd})$ is $\be_{ad}$ times 
the value of the product of $-h_{ac}\delta_{v_{bc}}$ and 
$\be_{cd}$ at 
the point $x(v_{da})\in\R$:  
\begin{equation*}
 \begin{split}
  P_{ad} m(-h_{ac}\delta_{v_{bc}},\be_{cd})
 & =\pm \left(e^{f_{ac}(x(v_{da}))-f_{ac}(x(v_{bc}))}\cdot 
 e^{f_{cd}(x(v_{da}))-f_{cd}(x(v_{cd}))}\right)\cdot \be_{ad} \\
 & =\pm \left(e^{-Y}\cdot e^{-Z}\right)\cdot\be_{ad},  
 \end{split}
\end{equation*}
where note that $f_{ac}(x(v_{bc}))-f_{ac}(x(v_{da}))=Y$ and 
$f_{cd}(x(v_{cd}))-f_{cd}(x(v_{da}))=Z$. 
Combining all these together, we obtain the first term in the right hand side of 
eq.(\ref{tree-transversal}):\  $\pm e^{-X-(Y+Z)}\be_{ad}$. 

In a similar way, one can see that the second term vanishes. 
The product $m(\be_{bc},\be_{cd})$ is proportional to $\delta_{v_{bc}}$, 
and its image by $h_{bd}$ is proportional to 
$e^{f_{bd}-f_{bd}(x(v_{bc}))}\cdot\vartheta^{-1}_{v_{bc}}$, whose value 
at the point $x(v_{da})\in\R$ is equal to zero. 
Therefore, $P_{ad} m(\be_{ab}, -h_{bd}\delta_{v_{bc}})$ vanishes 
due to the projection $P_{ad}$. 
\end{exmp}
This example shows how the transversal $A_\infty$ products are derived 
in general; the 
HPT machinery defines a higher product $m_k$ 
in terms of the sum of values 
associated to planar rooted $k$-trees over all the $k$-trees, 
but only the one compatible with the $k$-tree associated to the corresponding 
CC-polygon survives and produces the area of the CC-polygon. 
This phenomenon can be found in the original 
Morse homotopy theory \cite{fukayaAinfty,FO}, and also 
its extension \cite{KoSo} 
where the area of the polygon is taken into account as above.

The following is the first example of non-transversal products 
(\ref{str-const2}). 
\begin{exmp}\label{exmp:2}
$\vec{v}=(v_{ab},\vv_{bc},\vv_{cd},v_{da})$, 
$x(v_{bc})< x(v_{ab})< x(v_{da})< x(v_{cd})$. 
Let us calculate the non-transversal $A_\infty$-product 
$m_4(\be_{ab},\be_{bc},\be_{cd},\be_{da})$. This 
is again described as the sum of the values associated to 
trivalent planar rooted $4$-trees in the framework of HPT. 
In a similar way as in the transversal case above, 
the $4$-tree giving nonzero value is again the one corresponding to the 
CC-polygon only. 
Here, the $4$-tree corresponding to the CC-polygon $\vec{v}$ 
is obtained as follows. 
Connect the two degree zero points 
$v_{bc}$ and $v_{cd}$ with an interval. 
For each degree one point (in this case $v_{ab}$ and $v_{da}$), 
draw an interval, perpendicular to the $x$-axis, 
starting from the point and ending on the interval $(v_{bc}v_{cd})$. 
Then, we need a root edge; we add an edge, perpendicular to the $x$-axis, 
starting from a point on 
the interval $(v_{bc}v_{cd})$ between $x(v_{ab})$ and $x(v_{da})$ 
and ending on the interval $(v_{ab}v_{da})$ 
(Figure \ref{fig:non-transversal1}).  
\begin{figure}[h]
 \includegraphics{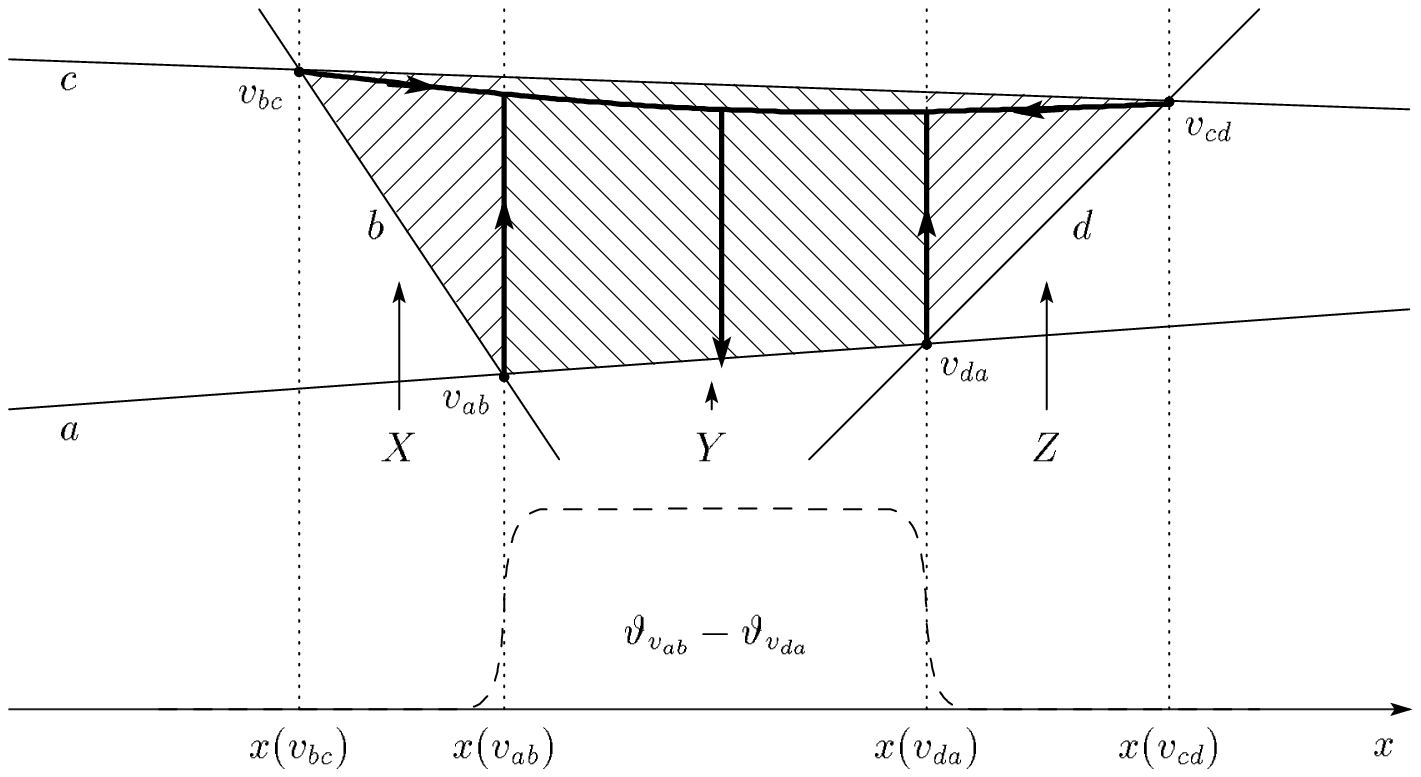}
 \caption{CC-polygon 
$\vec{v}=(v_{ab},\vv_{bc},\vv_{cd},v_{da})$. }\label{fig:non-transversal1}
\end{figure}
One can check that only the multilinear map 
corresponding to this $4$-tree is nonzero and it turns out to be 
$\pm e^{-(X+Y+Z)}(\vartheta_{v_{ab}}-\vartheta_{v_{da}})$. 
\end{exmp}
\begin{exmp}\label{exmp:3}
$\vec{v}=(\vv_{ab},v_{bc},\vv_{cd},v_{da})$, 
$x(v_{ab})<x(v_{bc})<x(v_{da})<x(v_{cd})$. Consider the non-transversal 
$A_\infty$-product 
$m_4(\be_{ab},\be_{bc},\be_{cd},\be_{da})$. In this case, 
there exist two choices of the $4$-trees corresponding to the 
CC-polygon $\vec{v}$. As in the previous example, 
we need to add an appropriate root edge. One can see that 
\begin{figure}[h]
 \includegraphics{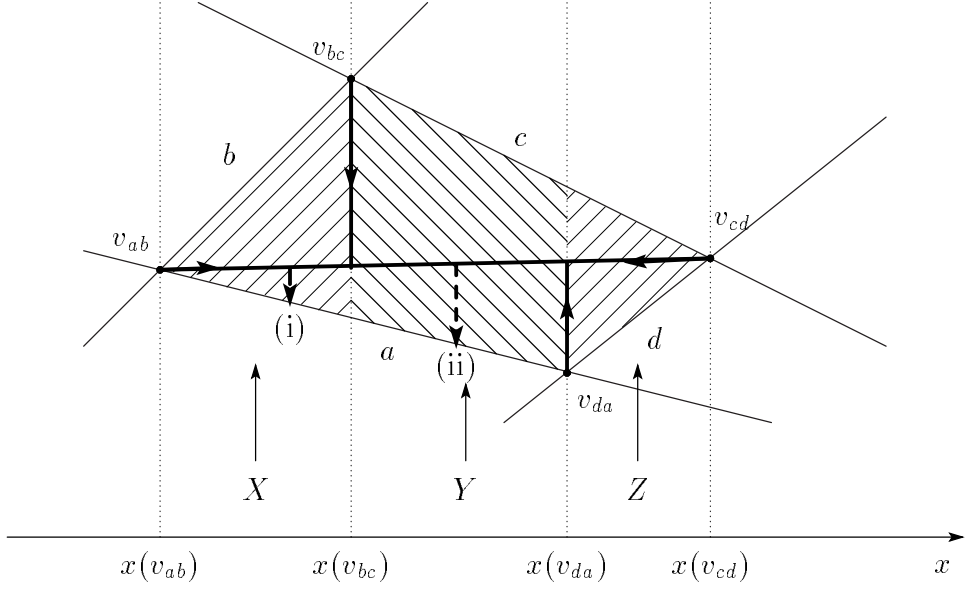}
 \caption{CC-polygon 
$\vec{v}=(\vv_{ab},v_{bc},\vv_{cd},v_{da})$ 
with two choices of the associated $4$-trees for the non-transversal 
$A_\infty$-product $m_4(\be_{ab},\be_{bc},\be_{cd},\be_{da})$. }
\label{fig:non-transversal2}
\end{figure}
there exist two choices (i) and (ii) of the root edge 
as in Figure \ref{fig:non-transversal2}. 
Actually, for the multilinear maps associated to $4$-trees by HPT, 
only those corresponding to these two $4$-trees 
give nonzero contribution. Recall that $P_{aa}=\Id_{\Omega'_{aa}}$ 
and one has 
\begin{equation*}
 \begin{split}
 m_4(\be_{ab},\be_{bc},\be_{cd},\be_{da})
  =& \pm m(\be_{ab},-h_{ba}m(\be_{bc},-h_{ca}m(\be_{cd},\be_{da})))\\
   & \pm m(-h_{ac}m(\be_{ab},\be_{bc}),-h_{ca}m(\be_{cd},\be_{da}))), 
 \end{split}
\end{equation*}
where the first (resp. second) term of the right hand side of the equation 
above corresponds to the $4$-tree with the root edge (i) (resp. (ii)). 
These two terms turn out to be 
\begin{equation*}
m(\be_{ab},-h_{ba}m(\be_{bc},-h_{ca}m(\be_{cd},\be_{da})))
  = \pm e^{-(X+(Y+Z))}(\vartheta_{v_{ab}}-\vartheta_{v_{bc}}), 
\end{equation*}
\begin{equation*}
m(-h_{ac}m(\be_{ab},\be_{bc}),-h_{ca}m(\be_{cd},\be_{da})))
  = \pm e^{-(X+Y+Z)}(\vartheta_{v_{bc}}-\vartheta_{v_{da}}), 
\end{equation*}
where the signs $\pm$ actually agree with each other and the result is 
\begin{equation*}
  m_4(\be_{ab},\be_{bc},\be_{cd},\be_{da})=
   \pm e^{-(X+Y+Z)}(\vartheta_{v_{ab}}-\vartheta_{v_{da}}). 
\end{equation*}
\end{exmp}
\begin{exmp}\label{exmp:4}
$\vec{v}=(\vv_{ab},(v_{bc})^{\otimes d},\vv_{cd},v_{da})$, 
$x(v_{ab})<x(v_{bc})<x(v_{da})<x(v_{cd})$. 
Consider the non-transversal $A_\infty$-product 
\begin{equation*}
m_{1+d+1}(\be_{ab},(\be_{bc})^{\otimes d},\be_{cd}):= 
m_{1+d+1}(\be_{ab},
\overbrace{\delta_{v_{bc}},\dots,\delta_{v_{bc}},\be_{bc},
\delta_{v_{bc}},\dots,\delta_{v_{bc}}}^d,
\be_{cd}). 
\end{equation*}
In fact, the result is independent 
of the order of $\delta_{v_{bc}}$'s and $\be_{bc}$, 
though $\delta_{v_{bc}}$'s in the left (resp. right) hand side of $\be_{bc}$ 
are elements in $V_{bb}^1$ (resp. $V_{cc}^1$). 
The corresponding CC-polygon is $\vec{v}$. Then, the situation is the same as 
Example \ref{exmp:1} for a transversal $A_\infty$-product 
except that we have $d$ elements associated to the point $v_{bc}$. 
One obtains 
\begin{equation*}
 m_{1+d+1}(\be_{ab},(\be_{bc})^{\otimes d},\be_{cd})
 =\pm P_{ad} m(-h_{ac}(w),\be_{cd}), 
\end{equation*}
where $h_{ac}(w)\in V_{ac}^0$ is given by 
\begin{equation*}
\pm h_{ac}m(\cdots -h_{ac}m(-h_{ac}m(-h_{ab} m(\cdots -h_{ab}m(\be_{ab},
 \delta_{v_{bc}}),\dots,\delta_{v_{bc}}),\be_{bc}),\delta_{v_{bc}}),
\dots,\delta_{v_{bc}}).
\end{equation*}
Since the final result is independent 
of the order of $\delta_{v_{bc}}$ and $\be_{bc}$, 
let us try to calculate this $w$ 
in the case all $\delta_{v_{bc}}$ is in the left hand side 
of the element $\be_{bc}$. Then, 
\begin{equation*}
 \begin{split}
 w=&  \pm -h_{ac} m( -h_{ab} m(\cdots -h_{ab}m(\be_{ab},
 \delta_{v_{bc}}),\dots,\delta_{v_{bc}}),\be_{bc}) \\
 & =\pm e^{-X}\cdot 
 h_{ac}m( \overbrace{(h_{ab}\delta_{v_{bc}}\cdot\cdots(h_{ab}\delta_{v_{bc}}\cdot
(h_{ab}\delta_{v_{bc}}}^{d-1} ))\cdots),\be_{bc}) \\
 & =\pm e^{-X} \ov{d!}(\vartheta_{v_{bc}})^d, 
 \end{split}
\end{equation*}
where, in the second line we omit denoting 
the product $m: {\Omega'}_{ab}^0\otimes {\Omega'}_{bb}^1\to {\Omega'}_{ab}^1$. 
In the third equality, 
recall that $\delta_{v_{bc}}=d\vartheta_{v_{bc}}$ and 
$h_{ab}d_{ab}(\vartheta_{v_{bc}})^k=(\vartheta_{v_{bc}})^k$, and then 
we used the formula $\delta_{v_{bc}}\cdot(\vartheta_{v_{bc}})^{k-1}
=(1/k)\, d(\vartheta_{v_{bc}})^k$ for $k=1,2,\dots$. 
The remaining calculation is the same as Example \ref{exmp:1} 
and we obtain $m_{1+d+1}(\be_{ab},(\be_{bc})^{\otimes d},\be_{cd})
=\pm (1/d!)\, e^{-(X+Y+Z)}\be_{ad}$. 
\end{exmp}
\begin{exmp}\label{exmp:5}
$\vec{v}=((\vv_{ab})^{\otimes (d_-,d_+)},v_{bc},\vv_{cd},v_{da})$, 
$x(v_{ab})<x(v_{bc})<x(v_{da})<x(v_{cd})$. 
Consider the non-transversal $A_\infty$-product 
\begin{equation*}
m_{d+1+1}((\be_{ab})^{\otimes (d_-,d_+)},\be_{bc},\be_{cd}):= 
m_{d+1+1}(\overbrace{\delta_{v_{ab}},\dots,\delta_{v_{ab}}}^{d_-},\be_{ab},
\overbrace{\delta_{v_{ab}},\dots,\delta_{v_{ab}}}^{d_+},\be_{bc},\be_{cd}), 
\end{equation*}
where $d:=d_-+d_+$. 
By HPT, there exist $d!/((d_-)!(d_+)!)$ number 
of $(d+1+1)$-trees which give nonzero contribution 
to the corresponding $(d+1+1)$-linear maps. 
The number $d!/((d_-)!(d_+)!)$ comes from the number of the orders that 
one acts $-h_{ab}m(\delta_{v_{ab}}, *)$ to $\be_{ab}$ $d_-$ times and 
$-h_{ab}m(*, \delta_{v_{ab}})$ to $\be_{ab}$ $d_+$ times. 
In fact, the result does not depend on the order and one can describe 
\begin{equation*}
 m_{d+1+1}((\be_{ab})^{\otimes (d_-,d_+)},(\be_{bc}),\be_{cd})
 =\frac{d!}{(d_-)!(d_+)!} 
 P_{ad} m(-h_{ac}m(w,\be_{bc}),\be_{cd}), 
\end{equation*}
where $w\in V_{ab}^0$ is given by 
\begin{equation*}
 \begin{split}
 & \pm \left(-h_{ab}m(\delta_{v_{ab}}, *) \right)^{d_-}\circ
 \left( -h_{ab}m(* , \delta_{v_{ab}})\right)^{d_+} \be_{ab} \\
 & =\pm \overbrace{(h_{ab}\delta_{v_{ab}}\cdot\cdots(h_{ab}\delta_{v_{ab}}
 \cdot(h_{ab}\delta_{v_{ab}}}^d))\cdots)\\
 & =\pm \ov{d!}\left(\vartheta_{v_{ab}}-\ov{2}\right)^d. 
 \end{split}
\end{equation*}
Here, in a similar way as in the previous example (Example \ref{exmp:4}), 
in the second line we omit denoting the product 
$m:{\Omega'}_{aa}^1\otimes {\Omega'}_{ab}^0\to {\Omega'}_{ab}^1$ 
or $m: {\Omega'}_{ab}^0\otimes {\Omega'}_{bb}^1\to {\Omega'}_{ab}^1$. 
The difference of the situation here from that 
in the previous example is that 
here $h_{ab}\delta_{v_{ab}}=\vartheta_{v_{ab}}-(1/2)$ 
and we have the formula 
$h_{ab}\delta_{v_{ab}}(\vartheta_{v_{ab}}-(1/2))^k
=(1/(k+1))(\vartheta_{v_{ab}}-(1/2))^{k+1}$, 
which is used in the third equality of the equation above. 
The remaining calculations are similar to those in Example \ref{exmp:1}; 
we finally obtain 
$P_{ad} m(-h_{ac}m(w,\be_{bc}),\be_{cd})
=\pm (1/(2^d d!))e^{-(X+Y+Z)}\be_{ad}$ and then 
\begin{equation*}
 m_{d+1+1}((\be_{ab})^{\otimes (d_-,d_+)},\be_{bc},\be_{cd})
 =\ov{2^d (d_-)!(d_+)!}\be_{ad}.  
\end{equation*}
\end{exmp}
\begin{exmp}\label{exmp:6} 
$\vec{v}=((\vv_{ab})^{\otimes(d_-,d_+)},
v_{bc},\vv_{cd},(v_{da})^{\otimes d'})$, 
$x(v_{bc})<x(v_{ab})<x(v_{cd})<x(v_{da})$. 
Consider the non-transversal $A_\infty$-product 
\begin{equation*}
m_{d+1+1+d'}((\be_{ab})^{\otimes (d_-,d_+)},
\be_{bc},\be_{cd},(\be_{da})^{\otimes d'}). 
\end{equation*}
This can be calculated by combining the arguments 
in Examples \ref{exmp:3}, \ref{exmp:4}, and \ref{exmp:5}. 
As in Example \ref{exmp:3}, the $A_\infty$-product is given as 
the sum 
\begin{equation}\label{exmp6-1}
 \begin{split}
 &  m_{d+1+1+d'}((\be_{ab})^{\otimes (d_-,d_+)},
\be_{bc},\be_{cd},(\be_{da})^{\otimes d'}) \\
 & = \pm m(w,-h_{ab}m(\be_{bc},w')) \pm m(-h_{ac}m(w,\be_{bc}),w'), 
 \end{split}
\end{equation}
where $w$ is just the $w$ in the previous example (Example \ref{exmp:5}), 
and $w'$ is given in a similar way as in Example \ref{exmp:4}: 
\begin{equation*}
  w= \pm\ov{(d_-)!(d_+)!}\left( \vartheta_{v_{ab}}-\ov{2}\right)^d,   \qquad 
  w'= \pm \ov{(d')!}(\vartheta^{-1}_{v_{da}})^{d'}. 
\end{equation*}
Here recall that $\vartheta^{-1}_{v_{da}}:=1-\vartheta_{v_{da}}$. 
Then, the two terms in the second line of eq.(\ref{exmp6-1}) turn out to be 
\begin{equation*}
 \begin{split}
 & m(w,-h_{ab}m(\be_{bc},w'))=
 \pm \ov{(d_-)!(d_+)!(d')!} 
 e^{-(X+Y+Z)}\ (\vartheta^{+1}_{v_{ab}})^d\cdot (\vartheta_{v_{bc}}^{-1}), \\
 & m(-h_{ac}m(w,\be_{bc}),w')= 
  \pm \ov{(d_-)!(d_+)!(d')!} 
 e^{-(X+Y+Z)} (\vartheta_{v_{bc}})\cdot (\vartheta^{-1}_{v_{da}})^{d'}. 
 \end{split}
\end{equation*}
The signs in fact agree with each other so that 
$(\vartheta_{v_{ab}})^d\cdot (\vartheta_{v_{bc}}^{-1})+
(\vartheta_{v_{bc}})\cdot (\vartheta^{-1}_{v_{da}})^{d'}
=\left((\vartheta_{v_{ab}})^d-(\vartheta_{v_{bc}})\right)+
\left((\vartheta_{v_{bc}})-\left(1-(\vartheta^{-1}_{v_{da}})^{d'}\right)\right)
= (\vartheta_{v_{ab}})^d\cdot (\vartheta_{v_{da}}^{-1})^{d'}$, 
and we finally obtain 
\begin{equation*}
 m_{d+1+1+d'}((\be_{ab})^{\otimes (d_-,d_+)},\be_{bc},\be_{cd},(\be_{da})^{\otimes d'})
 = \pm \ov{(d_-)!(d_+)!(d')!}
 (\vartheta^{+1}_{v_{ab}})^d\cdot (\vartheta_{v_{da}}^{-1})^{d'}. 
\end{equation*}
One can also check the case that $v_{bc}$ has 
multiplicity $d''$, $(v_{bc})^{\otimes d''}$, 
where, after using the following identity 
$\sum_{k=0}^dC_{k,d-k}(\vartheta_{v_{bc}})^k(\vartheta_{v_{bc}}^{-1})^{d-k}
=(\vartheta_{v_{bc}}+\vartheta_{v_{bc}}^{-1})^d=1$,  
we finally obtain just $(1/(d'')!)$ times the result above. 
\end{exmp}

\vspace*{0.3cm}

\noindent
$\bullet$\ {\bf The sign}\qquad 
The sign is determined precisely as follows. 
In order to simplify the sign in HPT, 
one may first consider the suspension $s(\cC_{DR}'(\fF_N))$ 
of $\cC_{DR}'(\fF_N)$. 
Then, apply HPT to $s(\cC_{DR}'(\fF_N))$ and obtain the $A_\infty$-products 
of $s(\cC(\cF_N))$. Finally, as the desuspension of 
$s(\cC(\cF_N))$ one obtains 
the $A_\infty$-products of $\cC(\fF_N)$. 

To see how the sign is determined, 
it is enough to demonstrate the calculations 
in the examples of transversal $A_\infty$-products below. 
As we saw in eq.(\ref{tree-transversal}), a transversal $A_\infty$-product $m_n$, 
$n\ge 2$, is 
described in terms of trivalent planar tree graphs, where 
the number of $m$ and $-h_{ab}$ for some $a\ne b\in\fF_N$ are 
$(n-1)$ and $(n-2)$, respectively. 
Since the sign problem for $m_2$ is obvious, 
let us consider the case $n\ge 3$. 
We obtain the sign from the following three parts. 
\begin{itemize}
 \item[$\circ$] 
For any product $m(w,w')$ in a tree graph, 
the degree $(|w|,|w'|)$ in $\cC'_{DR}(\fF_N)$ is $(0,1)$ or $(1,0)$. 
The suspension $s:\cC_{DR}'(\fF_N)\to s(\cC_{DR}'(\fF_N))$ 
leads to sign $(-1)^I$ for $I$ the number of the products $m(w,w')$ 
with degree of type $(0,1)$ (see eq.(\ref{suspension-sign})).

 \item[$\circ$] Associated to each internal edge, 
we have $-h_{ab}(\delta_{v})$ for some $a\ne b\in\fF_N$ 
and $v\in S_{all}$. 
The arrow of the internal edge is oriented from the left to the right or 
from the right to the left. 
Then, we have sign $(-1)^J$ where $J$ is the number of the internal edges 
oriented from the left to the right. 
(Compare this argument with eq.(\ref{h-ab}) and eq.(\ref{h-ba}) with $n=1$. )

 \item[$\circ$] In the process of the desuspension $s(\cC(\fF_N))\to\cC(\fF_N)$, 
an $A_\infty$-product $m_n(w_1,\dots,w_n)$ gets sign $(-1)^K$ with 
$K=n-i$ if $w_i$, for some $1\le i\le n$, is the only degree zero element, 
and $K=(n-i)+(n-j)$ 
if $w_i$ and $w_j$, for some $1\le i<j\le n$, 
are the only degree zero elements 
(see eq.(\ref{suspension-sign})). 
Note that by degree counting (Lemma \ref{lem:degree}) 
there are not more than two degree zero 
elements in $\{w_1,\dots,w_n\}$. 
\end{itemize}
Thus, $(-1)^{I+J+K}$ is the sign we finally obtain.

\begin{figure}[h]
 \begin{minipage}[c]{75mm}{
\begin{center}
 \includegraphics{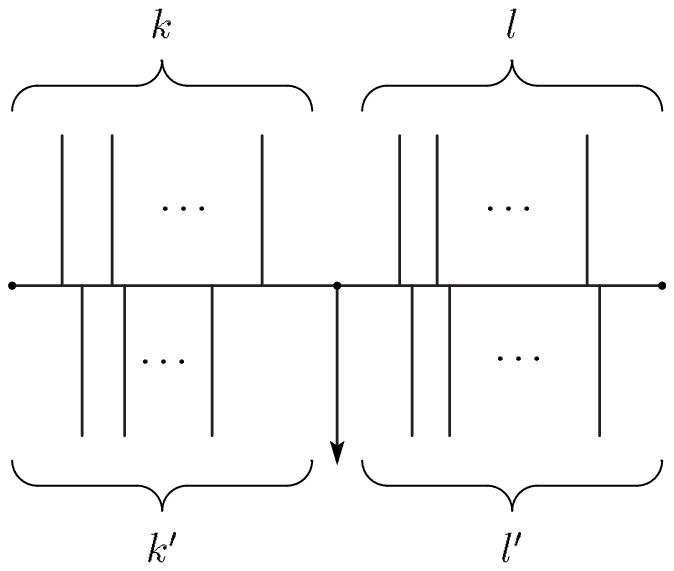}

 (a)
\end{center}}
 \end{minipage}
\quad 
 \begin{minipage}[c]{75mm}{
\begin{center}
 \includegraphics{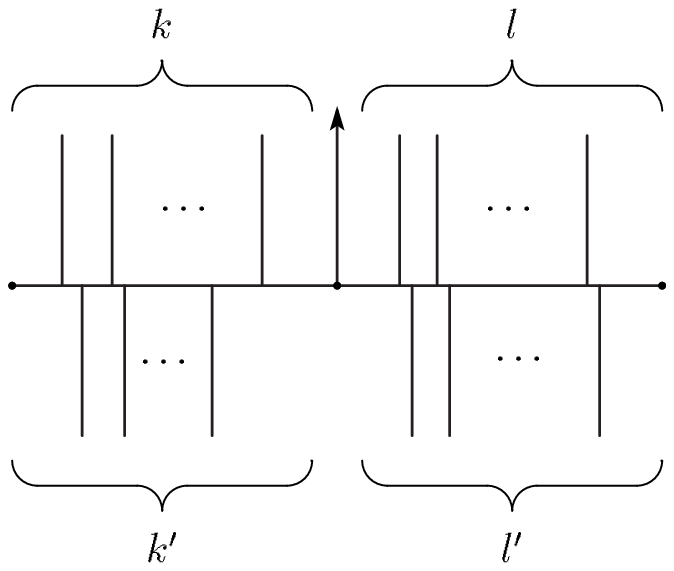}

 (b)
\end{center}}
 \end{minipage}
 \caption{The tree graphs corresponding to CC-polygons $\vec{v}$ 
with (a): $\sigma(\vec{v})=-1$ and (b): $\sigma(\vec{v})=+1$. }
\label{fig:sign}
\end{figure}

Let us consider the examples of the transversal $A_\infty$-products 
$m_n(w_1,\dots, w_n)$ with two degree zero elements $w_i$ and $w_j$ 
for some $1\le i<j\le n$. 
The corresponding tree graph is described 
as in Figure \ref{fig:sign} (a) and (b) 
when the corresponding CC-polygon $\vec{v}$ has $\sigma(v)=-1$ and 
$\sigma(v)=-1$, respectively. 
Here, note that $n=k+k'+l+l'+2$. 
($i=k'+1$ and $j=k'+1+k+l+1$ for case (a), and 
$i=l+1$ and $j=l+1+l'+k'+1$ for case (b). )
For the case (a), one obtains 
\begin{equation*}
 I=k+l'+1,\qquad J=k+k',\qquad K=l'+(k+l+1+l')
\end{equation*}
and hence the sign is $(-1)^{I+J+K}=(-1)^{k+k'+l+l'+2}=(-1)^n$. 
For the case (b), one obtains 
\begin{equation*}
 I=k+l'+1,\qquad J=k+k',\qquad K=k+(l'+k'+1+k)
\end{equation*}
and hence the sign is $(-1)^{I+J+K}=1$. 
One can see that in both cases the results agree with the definition of 
the transversal $A_\infty$-products in eq.(\ref{str-const1}). 
The calculation for the transversal $A_\infty$-product 
$m_n(w_1,\dots, w_n)$ with only one degree zero element 
in $\{w_1,\dots,w_n\}$ is similar.

 \section{Concluding remarks}
\label{sec:conclude}

We can apply the arguments in the present paper to two-tori directly. 
Then, we can discuss the homological mirror symmetry for two-tori 
including non-transversal $A_\infty$-products. 
In this case, note that we can include the identity morphism in the two-torus 
analog of the graded vector space $V_{aa}$, and further the DG category has 
a canonical nondegenerate inner product 
(canonical pairing of the Serre functor) 
defining cyclicity (see \cite{PoSc, hk:nchms} 
in noncommutative tori setting). 
Then, the dual of the identity morphism will be 
a natural representative of the cohomology of $V_{aa}^1$. 
Thus, if we start from finitely many objects, 
we can obtain an example 
of finite dimensional minimal $A_\infty$-algebras 
by applying HPT again to the two-torus analog 
of the graded vector space $V_{aa}$ 
of the $A_\infty$-category $\cC(\fF_N)$. 
From this viewpoint, 
the $A_\infty$-category $\cC(\fF_N)$ we constructed 
in this paper is an intermediate step. 
In particular, in $\R^2$ case, the $A_\infty$-structure is not equipped with 
cyclicity in the sense as in \cite{hk:nchms}; 
though cyclicity for an $A_\infty$-structure 
is defined by a {\em non-degenerate} inner product, 
the inner product defined naturally in this case becomes 
degenerate on $V_{aa}$. 

The generalization of the story of the present paper ($\R$ case) 
to $\R^n$ case is also an important issue, where, 
though the generalization of the DG-category $\cC_{DR}(\fF_N)$ 
is straightforward (see \cite{hk:ncDG,hk:nctheta}), 
we need to define a higher dimensional analog 
of the DG-category $\cC'_{DR}(\fF_N)$ so that 
HPT can be applied to it. 
The construction of the higher dimensional analogue of $\cC'_{DR}(\fF_N)$ 
is not straightforward, 
\footnote{For instance, even if we consider affine Lagrangians only, 
the orbits of the gradient flow are not affine. Since 
the action of the homotopy operator $h_{ab}$ is defined by the orbits, 
those non-affine orbits cause various subtleties 
as pointed out by K.~Fukaya to the author. 
Another approach to construct an $A_\infty$-structure of a Fukaya category 
on a torus fibration is discussed in \cite{fukaya:asymptotic} 
where we can avoid this kind of subtleties. 
}
but it seems not still impossible. 
This higher dimensional generalization also enables us to consider 
nontrivial noncommutative deformation of the $A_\infty$-categories 
(see \cite{hk:ncDG,hk:nctheta}). 

The reader might notice that elements in $V_{aa}^1$ 
played a special role in the present paper. 
In fact, the elements in $V_{aa}^1$ is related 
to open string background, \ie, 
the solutions of the Maurer-Cartan equation of the $A_\infty$-structure 
(see \cite{FOOO, hk:thesis}). 
In $\R^2$ case, the Maurer-Cartan equation will be trivial, 
which implies that 
all elements in $V_{aa}^1$ can be the solution 
of the Maurer-Cartan equation. 
Then, nontrivial deformation of lines to curves in $\R^2$ can also 
be taken into account in this framework. 

Finally, instead of the application to tori, 
we hope to apply the arguments of this paper to more general manifolds 
since the $A_\infty$-categories in this $\R^2$ case 
and their higher dimensional generalization, 
if it could be done, might be thought 
of a local construction of the $A_\infty$-categories 
which should be defined on the whole manifolds.


\end{document}